\documentclass{amsart}

\usepackage[all]{xy}

\newtheorem{theorem}{Theorem}[section]
\newtheorem*{maintheorem}{Theorem}
\newtheorem{lemma}[theorem]{Lemma}
\newtheorem{proposition}[theorem]{Proposition}
\newtheorem{corollary}[theorem]{Corollary}

\theoremstyle{definition}

\newtheorem{remark}[theorem]{Remark}
\newtheorem*{acknowledgement}{Acknowledgement}

\theoremstyle{remark}

\DeclareFontFamily{U}{wncy}{}
\DeclareFontShape{U}{wncy}{m}{n}{<->wncyr10}{}
\DeclareSymbolFont{mcy}{U}{wncy}{m}{n}
\DeclareMathSymbol{\Sh}{\mathord}{mcy}{"58}

\newcommand\mynote[1]{\marginpar{\ \\ \small \tt #1}}
\newcommand\bel[1]{{\mynote{#1}}\begin{equation}\label{#1}}
\newcommand\mylabel[1]{\label{#1}}

\newcommand{\ZZ}{\mathbb{Z}}
\newcommand{\QQ}{\mathbb{Q}}

\newcommand{\FF}{\mathbb{F}}

\newcommand{\PP}{\mathbb{P}}

\newcommand{\GG}{\mathbb{G}}

\newcommand  {\shA}     {\mathcal{A}}

\newcommand  {\shF}     {\mathcal{F}}

\newcommand  {\shI}     {\mathcal{I}}

\newcommand  {\shL}     {\mathcal{L}}

\newcommand  {\aid}     {\mathfrak{a}}

\renewcommand{\cong}    {\equiv}

\newcommand  {\dra}     {\dashrightarrow}
\newcommand  {\depth}   {\operatorname{depth}}

\newcommand  {\Div}     {\operatorname{Div}}

\newcommand  {\edim}    {\operatorname{edim}}
\newcommand  {\End}     {\operatorname{End}}

\newcommand  {\Ext}     {\operatorname{Ext}}

\newcommand  {\Fr}      {\operatorname{Fr}}

\newcommand  {\GL}      {\operatorname{GL}}

\newcommand  {\Hom}     {\operatorname{Hom}}
\newcommand  {\Hilb}    {\operatorname{Hilb}}

\newcommand  {\id}      {\operatorname{id}}

\newcommand  {\Km }     {\operatorname{Km}}

\newcommand  {\loc}     {{\operatorname{loc}}}

\newcommand  {\invlim}  {\varprojlim}

\newcommand  {\lra}     {\longrightarrow}

\newcommand  {\maxid}   {\mathfrak{m}}

\newcommand  {\nor}     {{\operatorname{nor}}}

\renewcommand{\O}       {\mathcal{O}}

\newcommand  {\Pic}     {\operatorname{Pic}}

\newcommand  {\sign}      {\left\{\pm 1\right\}}

\newcommand  {\Proj}    {\operatorname{Proj}}

\newcommand  {\quadand} {\quad\text{and}\quad}

\newcommand  {\ra}      {\rightarrow}

\newcommand  {\red}     {{\operatorname{red}}}

\renewcommand{\setminus}{\smallsetminus}

\newcommand  {\Sing}    {\operatorname{Sing}}

\newcommand  {\sm}      {\setminus}

\newcommand  {\Spec}    {\operatorname{Spec}}

\newcommand  {\Sym}     {\operatorname{Sym}}

\newcommand  {\TS}      {\operatorname{TS}}

\newcommand\kreis[1]{\ensuremath{\mathbin{\settowidth{\dimen7}{\mbox{$ \bigcirc$}}%
\makebox[0pt][l]{$\bigcirc$}\makebox[\dimen7]{#1}}}}

\def\mydate{\number\day\space\ifcase\month \or January\or February\or March\or 
April\or May\or June\or July\or
August\or September\or October\or November\or December\fi \space\number\year}

\DeclareFontFamily{U}{wncy}{}
\DeclareFontShape{U}{wncy}{m}{n}{<->wncyr10}{}
\DeclareSymbolFont{mcy}{U}{wncy}{m}{n}
\DeclareMathSymbol{\Sh}{\mathord}{mcy}{"58}

\usepackage{amscd,amssymb,amsmath,graphics}

\begin{document}

\title[Supersingular abelian surfaces]
      {The Hilbert scheme of   points  
       for supersingular abelian surfaces}

\author[Stefan Schroer]{Stefan Schr\"oer}
\address{Mathematisches Institut, Heinrich-Heine-Universit\"at,
40225 D\"usseldorf, Germany}
\curraddr{}
\email{schroeer@math.uni-duesseldorf.de}


\dedicatory{8 June 2006}

\begin{abstract}
We analyse the geometry of Hilbert schemes of points on abelian surfaces
and Beauville's generalized Kummer varieties in positive characteristics.
The main result is that, in characteristic two, the addition map from the Hilbert scheme of two points to
the abelian surface is a quasifibration, such that all
  fibers are nonsmooth. In particular, the corresponding generalized Kummer surface is 
nonsmooth, and  minimally elliptic singularities    occur in the supersingular case.
We unravel the structure of the singularities in dependence of $p$-rank and $a$-number
of the abelian surface.
To do so, we establish a McKay Correspondence for  Artin's wild involutions.
\end{abstract}

\maketitle
\tableofcontents

\section*{Introduction}
This paper circles, in positive characteristics, about the following subjects:
The \emph{McKay Correspondence}, \emph{Artin's Wild Involutions}, and the \emph{Hilbert--Chow morphism}.
My point of departure is Beauville's generalized Kummer construction, which works
as follows:

Fix a complex abelian surface $A$   and let $\Hilb^n(A)$ be its Hilbert scheme
of $n$ points. One knows that $\Hilb^n(A)$ is smooth with
trivial dualizing sheaf. It  is a crepant resolution of the symmetric product,
given by the Hilbert--Chow morphism $\gamma:\Hilb^n(A)\ra\Sym^n(A)$.
From this one gets an addition map $\Hilb^n(A)\ra A$, and
Beauville \cite{Beauville 1983} introduced the \emph{generalized Kummer variety} $\Km^n(A)$ as the    fiber 
of the addition map   over the origin.
It turns out that $\Km^n(A)$ is   smooth, and its dualizing sheaf is trivial as well.
In fact, generalized Kummer varieties are important examples of   \emph{hyperk\"ahler manifolds}.

The same construction works over ground fields of characteristic $p>0$.
The Hilbert scheme $\Hilb^n(A)$ is still smooth with trivial dualizing sheaf.
For the generalized Kummer varieties, however, entirely new geometric phenomena
arise: As we shall see, $\Km^n(A)$ is not necessarily smooth, and may even be nonnormal.
The goal of this paper is to study the addition map  
and the generalized Kummer variety  
in the simplest accessible case, namely in characteristic $p=2$
for $n=2$ points. The first main result is the following:

\begin{maintheorem}
In characteristic two, all fibers of the addition map $\Hilb^2(A)\ra A$
are nonsmooth. They are always geometrically reduced, and geometrically normal if and only if 
the abelian variety $A$ is not superspecial.
\end{maintheorem}

Whence $\Hilb^2(A)\ra A$ is an example of a \emph{quasifibration},
that is, its generic fiber is regular but not geometrically regular.
Such a startling violation of Sard's Lemma is only possible in  positive characteristics.
So far, there is little systematic study of quasifibrations, except for the special case of
quasielliptic  surfaces, which  play an important role
in the Enriques classification in positive characteristics \cite{Bombieri; Mumford 1976}.

The generalizsed Kummer surface  $\Km^2(A)$  is    related to 
the classical Kummer surface, which is the quotient  $A/\left\{\pm 1\right\}$ of the abelian surface
by the sign involution. For such quotients Artin's classification \cite{Artin 1975} of 
involutions on surfaces in characteristic two applies.
Shioda \cite{Shioda 1974} and Katsura \cite{Katsura 1978}  proved that  the singularities on the normal surface $A/\left\{\pm 1\right\}$
are certain rational or elliptic double points.
This is in startling contrast to the complex situation, where we always have sixteen ordinary double points.
The second main result of this paper  is a description of the singularities on $\Km^2(X)$
in relation to the singularities on $A/\sign$:

\begin{maintheorem}
Suppose the ground field $k$ is of characteristic $p=2$. Then
Beauville's Kummer surface $\Km^2(A)$ is crepant partial  desingularization of the
classical Kummer surface $A/\left\{\pm 1\right\}$ obtained by blowing-up
the the schematic image of the fixed scheme on $A$.
\end{maintheorem}

The precise structure of the singularities on $\Km^2(A)$ will be determined in 
Sections \ref{the sign involution}, \ref{minimally elliptic}, and \ref{superspecial case}.
The following table   gives a rough idea of the situation:
$$
\hspace{-2.8em}\begin{array}[t]{c|c|c|c|c}
\text{$p$-rank or $a$-number of $A$}         &  \sigma=2    &  \sigma=1  &  \sigma=0,a=1 &  a=2 \\
\hline&&&\\[-2ex]
\text{Singularities on $A/\sign$}  &  4D_4^1 & 2D_8^2 & \text{elliptic double point} & \text{elliptic double point}\\
\hline&&&\\[-2ex]
\text{Singularities on $\Km^2(A)$} &  12A_1  &  2A_3+2D_4^0  &  \text{elliptic triple point} & \text{nonnormal} \\
\end{array}
$$

\medskip
Here the first row contains the two basic numerical invariants of   abelian varieties in positive characteristic, namely the $p$-rank $\sigma$ and the $a$-number $a$. The upper indices in $D_n^r$ determines the isomorphism class of rational double points
of type $D_n$ in characteristic two, according to Artin's analysis \cite{Artin 1977}. 
The supersingular case is most challenging: Here our analysis
depends on Laufer's theory of \emph{minimally elliptic singularities} \cite{Laufer 1977}.

The existence of a crepant partial resolution   holds true in general for quotients of surfaces by involutions in characteristic two,
and is closely related to $G$-Hilbert schemes.
Recall   the complex \emph{McKay Correspondence} in dimension two
 was established in various degrees of generality:
Ito and Nakamura \cite{Ito; Nakamura 1999} showed that the minimal resolution of singularities for rational double points  is isomorphic to a suitable $G$-Hilbert scheme. 
This was extended by Kidoh \cite{Kidoh 2001} to 
cyclic quotients singularities, and by Ishii \cite{Ishii 2002} to arbitrary quotient singularities.

The situation appears to be more complex in positive characteristics.
Suppose that $S$ is a quasiprojective smooth surface in characteristic $p=2$, endowed with an action
of $G=\sign$ having a single fixed point $s\in S$. Let $T=S/G$ be the quotient surface.
Then the image $t=q(s)$ of the fixed point under the quotient map $q:S\ra T$
is an isolated singularity.
The third main result of this paper describes the McKay Correspondence in this situation:

\begin{maintheorem}
The blowing-up $g:T'\ra T$ of the image of the fixed scheme $q(S^G)\subset T$
is a crepant partial resolution with $R^1g_*(\O_{T'})=0$. The scheme $T'$ is 
isomorphic to the reduced $G$-Hilbert scheme $\Hilb^G_\red(S)$.
\end{maintheorem}

The scheme $T'$ is usually nonnormal, and the nonreduced $G$-Hilbert scheme $\Hilb^G(S)$
usually contains embedded components. 
Note that the two descriptions $T'=\Hilb^G_\red(S)$ are entirely different:
The first is effective for explicit computations, the second is useful for  
theoretical considerations.

\medskip
Here is an outline of the paper:
In Section \ref{wild involutions}, we recall Artin's structure result on
wild involutions, and introduce our crepant partial resolution $g:T'\ra T$.
In Section \ref{maps hilbert} we show that $T'$ is isomorphic to a closed subscheme of the Hilbert scheme $\Hilb^2(S)$,
which  depends on an   infinitesimal computation.
In Section \ref{g-hilbert schemes}, we identify $T'$ with the reduced $G$-Hilbert scheme
$\Hilb^G_\red(S)$. In Section \ref{example: rdp}, we discuss the case that the singularity $\O_t$
is a rational double point in detail.
After that, we turn to abelian surfaces: In Section \ref{the sign involution}, I apply the
preceding results to the sign involution on an abelian surface $A$. 
The case of abelian varieties with nonzero $p$-rank is comparatively straightforward.
The supersingular case is treated in the next two sections:
In Section \ref{minimally elliptic}, we use the theory of minimally elliptic singularities
to determine the singularities and the crepant partial resolution for $A/\sign$ 
in the case that $A$ is supersingular but not superspecial.
In Section \ref{superspecial case}, we unravel the superspecial case.
Here the crepant partial resolution is nonnormal. It is obtained from
the blowing-up of the reduced singular point, which is nonnormal as well,
by a kind of infinitesimal flip.
In Section \ref{serre conditions}, we analyse symmetric products and
the extent of their Cohen--Macaulayness. 
In Section \ref{symmetric products}, we consider symmetric products
of abelian varieties and the effect of the addition map on the geometry
of the symmetric product. 
In Section \ref{hilbert-chow}, we gather all our previous results,
and analyse the Hilbert-Chow morphism.

\begin{acknowledgement}
I wish to thank Torsten Ekedahl for helpful discussions.
\end{acknowledgement}

\section{Artin's wild involutions}
\mylabel{wild involutions}

Fix a ground field $k$ of characteristic $p=2$, and
let $S$ be a quasiprojective smooth surface, endowed with
an involution $\iota:S\ra S$.
In other words, the cyclic group  $G=\left\{\pm 1\right\}$ of order two
acts on $S$. Then the quotient $T=S/G$ is a quasiprojective normal surface.
Let $q:S\ra T$ be the quotient morphism.
To simplify, we assume that $G$ acts freely except for 
a single rational fixed point $s\in S$.
Let $t\in T$ be the image of this fixed point, such that 
$\Sing(T)=\left\{t\right\}$.

The goal of this section is to study the singularity $t\in T$ in terms of
the  blowing-up $T'\ra T$ with center the image of the fixed scheme $S^G\subset S$. 
We shall see that this blowing-up 
behaves cohomologically like the resolution of singularities for rational double points.
Therefore, it might be seen as a partial resolution of singularities.
It is, however,   usually   not normal.
In the next section, we shall identify our partial resolution with the underlying 
reduced subscheme of the $G$-Hilbert scheme of $S$.

To start with, I recall Artin's work on   involutions on surfaces in characteristic two.
In contrast to the case of complex numbers or odd characteristics,
the involution $\iota$ acting on the complete local ring $\O_{S,s}^\wedge$ 
is in general neither linearizable nor splits into a product action,
such that no obvious description of the quotient springs to mind.
However, Artin \cite{Artin 1975} obtained the following:

\begin{proposition}
\mylabel{normal form}
There is a parameter system $x,y\in\O_{S,s}^\wedge$,   a regular
system of parameters $u,v\in\O_{S,s}^\wedge$, and a parameter system
$a,b\in k[[x,y]]$ so that
\begin{equation}
\label{artin description}
u^2+au +x =0 \quadand
v^2+bv + y =0,
\end{equation}
and we have $\O_{T,t}^\wedge=k[[x,y,z]]/(z^2+abz+xb^2+ya^2)$.
\end{proposition}

Some comments might be helpful.
First note that (\ref{artin description}), together with the Implicit Function Theorem for formal power series
(\cite{A 4-7}, Corollary to Proposition 10 in Chapter IV, \S4, No.\ 7),
implies that $x,y$ are in a unique way formal power series in the indeterminates $u,v$, and hence may be regarded
as elements in $\O_{S,s}^\wedge$.
Next, we observe that the involution $\iota$ must interchange the roots of
the two quadratic equations in (\ref{artin description}), whence is given by
\begin{equation}
\label{artin action}
u\longmapsto u+a\quadand v\longmapsto v+b.
\end{equation}
Consequently, we have 
\begin{equation}
\label{artin variables}
x=u^2+au=u\bar{u}\quadand
y=v^2+bv=v\bar{v}\quadand
z=u\bar{v} + \bar{u}v=ub+va,
\end{equation}
 where   $\bar{v}=\iota(v)$ etc.\ denotes the
action of the involution.

What is the structure of the singularity $\O_{T,t}$ and its resolution of singularities $\tilde{T}\ra T$?
I find it difficult to   make   general statements.
However, the   exceptional divisor  $E\subset\tilde{T}$ has the following property:

\begin{proposition}
\mylabel{exceptional tree}
The Picard scheme $\Pic^0_{E/k}$ is unipotent.
\end{proposition}

\proof
We shall use the  \emph{local fundamental group} $\pi_1^\loc(\O_{T,t})$, which is by definition the
fundamental group of the pointed spectrum $\Spec(\O_{T,t})\setminus\left\{t\right\}$.
For the problem at hand,  we may assume that the ground field $k$ is separably closed, 
and replace $T$ by the formal completion at $t\in T$. Then $\pi_1^\loc(\O_{T,t})$ is cyclic of order two.
Seeking a contradiction, we now assume that $\Pic^0_{E/k}$ is not unipotent.
Then $\Pic(E)$ contains nonzero elements of finite order prime to $p=2$,
say of order three.
Such an  elements element extends to a an element of order three in $\invlim\Pic(nE)$, which
follows from the exact sequences
$$
H^1(E,\O_{E}(-nE))\lra\Pic((n+1)E)\lra\Pic(nE)\lra H^2(E,\O_E(-nE)).
$$
By Grothendieck's Existence Theorem, this corresponds to  an element
in $\Pic(\tilde{T})$ of order three.
In turn, we obtain an invertible sheaf $\shL$ on $U=T\setminus\left\{t\right\}$ of order three.
Choosing a trivialization $\shL^{\otimes 3}$, we endow
$\shA=\O_U\oplus\shL\oplus\shL^{\otimes 2}$ with an algebra structure, and in turn
have a connected finite \'etale covering of degree three, contradiction.
\qed

\medskip
It follows that the integral components $E_i\subset E$ have genus zero, for otherwise the Picard scheme
would contain abelian varieties. Moreover, the intersection graph for the $E_i$ must be a tree,
because otherwise the Picard scheme would contain copies of $\GG_m$. Compare
the discussion in \cite{Bosch; Luetkebohmert; Raynaud 1990}, Chapter 9.

It is easy to determine the schematic fiber $q^{-1}(t)\subset S$
of the singular point $t\in T$:

\begin{lemma}
\mylabel{fiber point}
The ideal of the schematic fiber $q^{-1}(t)\subset S$ is 
generated by the elements $u^2,v^2\in \O_{S,s}^\wedge$. 
\end{lemma}

\proof
We have to check the equality of ideals $(x,y,z)\O_{S,s}^\wedge=(u^2,v^2)$.
The inclusion ``$\supset$'' follows directly from (\ref{artin description}).
To check the reverse inclusion, we  use (\ref{artin description}) to infer 
\begin{equation*}
\label{artin developement}
x\cong u^2 + a_xu^3+a_yuv^2\quadand
y\cong v^2+b_xu^2v+b_yv^3 \text{\quad modulo $(u,v)^4$},
\end{equation*}
where $a_x,a_y,b_x,b_y\in k$ are the coefficients of the linear monomials in the development
$a=a_xx+a_yy+O(2)$ and $b=b_xx+b_yy+O(2)$.
From this, together with $z=ub+va$ from (\ref{artin variables}), we deduce $x,y,z\in (u^2,v^2)$.
\qed

\begin{remark}
We may describe the ideal of the fiber without passing to 
the formal completion as follows:
If $u',v'\in\O_{S,s}$ is any regular parameter system, then
we have $(u^2,v^2)=(u'^2,v'^2)$ inside the formal completion.
Hence $(u'^2,v'^2)\subset\O_{S,s}$ is the ideal of the fiber
$q^{-1}(t)\subset S$. This leads to a coordinate free description of
the ideal for $q^{-1}(t)\subset S$ as the \emph{Frobenius power}
$\maxid_{s}^{[2]} =(f^2\mid f\in\maxid_s)$.
\end{remark}

Let $S^G\subset S$ be the \emph{fixed scheme} of the  $G$-action.
This is the largest closed subscheme on which the $G$-action is trivial.
In light of (\ref{artin action}), its ideal is generated by the parameter system $a,b\in \O_{S,s}^\wedge$.
Note that the Artin scheme $S^G$ is never reduced.

\begin{lemma}
The schematic image $q(S^G)\subset T$ of the fixed scheme $S^G\subset S$ under the quotient map $q:S\ra T$ 
is defined by the primary ideal  $(a,b,z)\subset\O_{T,t}^\wedge$.
\end{lemma}

\proof
Let $\aid=\O_{T,t}^\wedge\cap(a,b)\O_{S,s}^\wedge$ be the ideal of  the schematic image $q(S^G)\subset T$.
We have $z=ub+va$ by (\ref{artin variables}), hence $z\in \aid$, and therefore
$(a,b,z)\subset \aid$.
To check the reverse inclusion, let $f(x,y)+zg(x,y)=r(u,v)a+s(u,v)b$ be an element from the ideal $\aid$.
Since we already know $z\in \aid$, we may as well assume $g=0$. It remains  to check that $f$
vanishes in $k[[x,y]]/(a,b)$. But this is true,
because $k[[x,y]]\subset k[[u,v]]$ is faithfully flat  and $f$ vanishes
in $k[[u,v]]/(a,b)$.
\qed

\medskip
It follows that the ideal $(a,b,z)\subset\O_{T,t}^\wedge$ has a coordinate-free description
as the ideal of the schematic image of the fixed scheme.
Let  $g:T'\ra T$ be the blowing-up of this ideal, or equivalent the blowing-up with center $q(S^G)\subset T$.  
The following result asserts that this morphism behaves cohomologically like
the  resolution of rational double points.
One should keep in mind, however, that the scheme $T'$ is usually not normal,
as we shall see below.

\begin{theorem}
\mylabel{first blowing-up}
The scheme $T'$ is locally of complete intersection.
The fiber $g^{-1}(t)$ is isomorphic to the infinitesimal extension
of $\PP^1$ by $\O_{\PP^1}(-1)$. Furthermore, we have $R^1g_*(\O_{T'})=0$, and the relative 
dualizing sheaf $\omega_{T'/T}$ is trivial.
\end{theorem}

\proof
To verify this we may assume   that our scheme $T$   actually equals the spectrum of the ring $k[x,y,z]/(z^2+abz+xb^2+ya^2)$
and forget about  formal completions, which simplifies notation a little bit.
The blowing-up $g:T'\ra T$ of the ideal $(a,b,z)\subset\O_{T,t}$ is covered by three affine charts: The $a$-chart,
the $b$-chart, and the $z$-chart.

The $z$-chart is generated by variables $x,y,z,a/z,b/z$,
subject to  the relations
$1+(a/z)(b/z)z+x(b/z)^2+y(a/z)^2=0$ and $a/z\cdot z=a$ and $b/z\cdot z=b$.
The exceptional divisor is given by the additional relation $z=0$, which is easily seen to be empty.
We may therefore concentrate on the $a$-chart, the situation for the $b$-chart being
symmetric.

The $a$-chart is generated by four variables $x,y,b/a,z/a$ modulo two relations
\begin{equation}
\label{a-chart}
(z/a)^2+a\cdot b/a\cdot z/a + x(b/a)^2 + y=0\quadand 
b/a\cdot a=b.
\end{equation}
These equations   clearly correspond  to a regular sequence. It follows 
that $T'$ is locally of complete intersection.

We next examine the fiber $F=g^{-1}(t)$, which is of codimension one.
Obviously, it is
covered by two affine charts: The $a$-chart has generators $b/a,z/a$, with only relation
$(z/a)^2=0$. The  $b$-chart has generators $a/b,z/b$ with relation $(z/b)^2=0$.
The infinitesimal generators are related by $(z/a)=(b/a)(z/b)$ on the overlap.
We infer that $F=g^{-1}(t)$ is an infinitesimal extension of the projective
line $\PP^1$ by the invertible sheaf $\O_{\PP^1}(-1)$.
It follows from \cite{Bayer; Eisenbud 1995}, Theorem 1.2 that such extensions are unique up to
isomorphism. We note in passing that $F$ must be isomorphic to a nonreduced
quadric in $\PP^2$.

To proceed, consider the Cartier divisor $C\subset T'$ whose ideal is the tautological
sheaf $\O_{T'}(1)\subset\O_{T'}$ attached to the blowing-up. Note that $C$ is an infinitesimal extension of the fiber
$F=g^{-1}(t)$. The $a$-chart for $C$ has
generators $x,y,b/a,z/a$, modulo the relations $a$, $b$, and
$(z/a)^2 + x(b/a)^2 + y$. From this we infer that $\O_C$ has a decomposition series
\begin{equation}
\label{filtration}
0=\shI_0\subset\shI_1\subset\ldots\subset\shI_l=\O_C
\end{equation}
with factors $\shI_i/\shI_{i-1}\simeq\O_F$. The length $l$ of the composition series
is also the length of the Artin ring $\O_{T,t}/(a,b,z)\simeq k[[x,y]]/(a,b)$, which defines the center for
our blowing-up $T'\ra T$. Now let $C_i\subset C$ be 
the closed subscheme defined by $\shI_i\subset\O_F$. 
We infer that an  invertible $\O_C$-module $\shL$ with $\shL\cdot F\geq 0$
has $H^1(C,\shL)=0$; this follows inductively from  the exact sequences
$$
H^1(F,\shL_F)\lra H^1(C_{i+1},\shL_{C_{i+1}})\lra H^1(C_i,\shL_{F_i}).
$$
In particular, we have $H^1(C,\O_C)=0$.

Next, let $C_n$ be the $n$-th infinitesimal neighborhood of
the  Cartier divisor $C=C_0$, where $\shI=\O_{T'}(1)$, and $F_n$ be the $n$-th infinitesimal neighborhood of the fiber $F=g^{-1}(t)$. Given $n\geq 0$, there is an $m\geq n$ with
 $C_n\subset F_m$ and $F_n\subset C_m$, and the restriction maps
$$
H^1(F_m,\O_{F_m})\lra H^1(C_n,\O_{C_n})\quadand
H^1(C_m,\O_{C_m})\lra H^1(F_n,\O_{F_n})
$$
are surjective. This implies that the two inverse systems of groups $H^1(F_m,\O_{F_m})$ and $H^1(C_m,\O_{C_m})$ have isomorphic inverse limits.
Using the Theorem on Formal Functions, we infer that the  canonical map
$$
R^1g_*(\O_{T'})\lra\invlim H^1(C_n,\O_{C_n})
$$
is bijective.
On the other hand, the short exact sequence of coherent sheaves supported by the exceptional locus
$0\ra \shI^{n}/\shI^{n+1}\ra\O_{F_{n+1}}\ra\O_{C_n}\ra 0$
yields an exact sequence
$$
H^1(C,\shI^n/\shI^{n+1})\lra H^1(C_{n+1},\O_{C_{n+1}})\lra H^1(C_n,\O_{C_n})\lra 0.
$$
The sheaf  $\shI^n/\shI^{n+1}\simeq\O_C(n)$ is invertible and ample,
and we saw in the preceding paragraph that this implies $H^1(C,\shI^n/\shI^{n+1})=0$.
Using induction on $n$, we infer that $H^1(C_n,\O_{C_n})$  vanishes.
Whence $R^1g_*(\O_{T'})=0$.

It remains to check that the relative dualizing sheaf is trivial.
Since it is trivial outside the fiber $F=g^{-1}(t)$, there is a Cartier divisor
$D\in\Div(T')$ supported by $F$ with $\omega_{T'/T}=\O_{T'}(D)$.
Our task is to show $D=0$.
Consider the Cartier divisor $C\subset T'$. Its dualizing sheaf is
$\omega_C=\O_C(C+D)$, by relative duality for the inclusion $C\subset T'$.
Serre duality gives $\deg(\omega_C)=-2\chi(\O_C)$.
Using the decomposition series (\ref{filtration}), we infer 
that $\chi(\O_C)=-l$.
Hence $\deg(\omega_C)=-2l$.
On the other hand, the invertible $\O_F$-module $\O_F(-C)$
is generated on the $a$-chart by $a$, and on the $b$-chart by $b$,
with $a=a/b\cdot b$ on the overlap.
It follows that $\deg\O_F(-C)=2$, and whence $\deg\O_C(C)=-2l$.
Consequently, $\O_C(D)$ had degree zero. The curve $C$ is irreducible, and  its Picard scheme has  tangent space $H^1(C,\O_C)=0$.
Whence we actually have  $\O_C(D)\simeq\O_C$, and the  same applies to all infinitesimal neighborhoods $C_n$.
Applying the Theorem of Formal Functions once more, we see that
$f_*\O_{T'}(D)$ is an invertible sheaf,  which  is obviously trivial outside $t\in T$.
Now \cite{Hartshorne 1994}, Theorem 1.12 tells us that $f_*\O_{T'}(D)$ is   trivial.
This implies $\O_{T'}(D)\simeq\O_{T'}$.
\qed

\medskip
A partial resolution whose relative dualizing sheaf is trivial is  called \emph{crepant}.
Although this terminology is usually applied   to normal partial resolutions,
we shall also say that our blowing-up $T'\ra T$ is crepant.

In characteristic zero, any quotient singularity is rational, according
to Hochster and Eagon \cite{Hochster; Eagon 1971}.
This does not hold true in positive characteristics. 
In case our quotient singularity  happens to be rational, this has the following
consequence:

\begin{corollary}
\mylabel{rational singularity}
Suppose the $\O_{T,t}$ is a rational singularity.
Then $\O_{T,t}$ is a rational double point, the surface
$T'$ is normal with only rational double points, and
the minimal resolution $\tilde{T}\ra T$ factors over our  partial resolution $T'\ra T$.
\end{corollary}

\proof
Recall that rational double points are precisely the rational Gorenstein singularities.
The singularity $\O_{T,t}$ is a complete intersection, and rational by assumption,
whence a rational double point.
Let $T^\flat\ra T'$ be the normalization, and $\tilde{T}^\flat\ra T^\flat$ be the minimal resolution of singularities.
Then we have a commutative diagram
$$
\begin{CD}
\tilde{T}^\flat @>>> T^\flat\\
@VVV @VVV\\
\tilde{T} @>>> T.
\end{CD}
$$
Suppose that the normalization $\nu:T^\flat\ra T'$ is not an isomorphism. Then it is not an isomorphism
over the exceptional divisor $C\subset T'$, because $T'$
satisfies Serre's Condition $(S_2)$. The relative dualizing sheaf $\omega_{T^\flat/T'}$
is given by the conductor ideal for the inclusion $\O_{T'}\subset\O_{T^\flat}$.
Using $K_{T'/T}=0$, we conclude that $K_{\tilde{T}^\flat/T}$ is not effective.
On the other hand, we have $K_{T'/T}=0$ since $\O_{T,t}$ is a rational double point, 
and $K_{\tilde{T}^\flat/T'}\geq 0$  because 
the morphism $\tilde{T}^\flat\ra T'$ decomposes into a sequence of blowing-ups of reduced points.
Hence $K_{\tilde{T}^\flat/T}\geq 0$, contradiction. Therefore $T^\flat=T'$ must be normal.

Finally, suppose that $\tilde{T}^\flat\ra\tilde{T}$ is not an isomorphism.
Then $K_{\tilde{T}^\flat/\tilde{T}}>0$. On the other hand, we have $K_{T^\flat/T}=0$ 
by Theorem \ref{first blowing-up},
and $K_{\tilde{T}^\flat/T^\flat}\leq 0$ because the resolution $\tilde{T}^\flat\ra T^\flat$ is minimal.
This gives again a contraction.
\qed

\medskip
We will examine the case of rational double points at length in Section \ref{example: rdp}.
As we shall see later, neither normality nor factorization holds true with minimally elliptic
singularities instead of rational singularities.

Let me now discuss  the question what Cartier divisors $D\subset Y'$ are supported on the
Weil divisor $F=g^{-1}(t)$. Given such a Cartier divisor, we have an equality 
$D=m F_\red$ of Weil divisors for some integer $m\geq 0$. The multiplicity $m$ is related to
the length $l\geq 1$ of the Artin algebra $\O_{T,t}^\wedge/(a,b,z)$,
which defines the center for the blowing-up $T'\ra T$, as follows:

\begin{corollary}
\mylabel{weil cartier}
The multiplicity $m$ of any    Cartier divisor $D\subset Y'$ supported by $F=g^{-1}(t)$ is a multiple
of   $2l$.
\end{corollary} 

\proof
Let $C\subset T$ be the Cartier divisor corresponding to the invertible sheaf $\O_{T'}(1)$
attached to the blowing-up $T'\ra T$. If follows from the computations in the preceding proof
that $F=2 F_\red$, and $C=lF$, and $F_\red\cdot \O_{Y'}(1)=-1$. 
We infer that $-m=mF_\red\cdot C=D\cdot C=D\cdot 2l F_\red$. The assertion follows,
because the number $D\cdot F_\red$ is an integer.
\qed

\medskip
As a consequence, there is no simple relationship between our blowing-up
of $q(S^G)\subset T$ and the blowing-up of the reduced singular point $t\in T$:

\begin{corollary}
\mylabel{partial factors}
Our partial resolution $g:T'\ra T$ factors over the blowing-up $T''\ra T$ of the
singularity $t\in T$ if and only if $(a,b,z)=(x,y,z)$. In this case, the two blowing-ups
coincide.
\end{corollary}

\proof
The condition is obviously sufficient.
Conversely, suppose there exist a factorization $T'\ra T''$.
The universal property of the blowing-up $T''\ra T$ implies that the Weil divisor
$F=g^{-1}(t)\subset T$ is Cartier. We have $F=2F_\red$, and the preceding
corollary tell us that $l=1$, whence $(a,b,z)=(x,y,z)$.
\qed

\medskip
The following observation will be important later:
Suppose our partial resolution  $T'$ is normal, and let $r:\tilde{T}\ra T'$ be 
the minimal resolution of singularities. For any Weil divisor $D$ on $T'$,
we then have the pullback $r^*(D)\in\Div(\tilde{T})\otimes\QQ$ in the sense of 
Mumford \cite{Mumford 1961}. We call $D$  \emph{numerically Cartier}
if the $\QQ$-divisor  $r^*(D)$ has integral coefficients. Then for any Weil divisor $D'$ on $Y'$,
the intersection number $D\cdot D'\in\QQ$ is   an integer as well.
In the special case that  $D=mF_\red$ is  supported
on the exceptional locus, the same proof as for Corollary \ref{weil cartier} gives:

\begin{corollary}
\mylabel{numerically cartier}
Suppose that the scheme $T'$ is normal, and that the Weil divisor $mF_\red$ is numerically Cartier.
Then the integer $m$ is a multiple of $2l$.
\end{corollary}

To finish this section, I want to clarify the dependence of the singular locus $\Sing(T')\subset T'$ 
on the parameters $a,b\in k[[x,y]]$.
Let $a_x,a_y,b_x,b_y\in k$ be the scalars describing the linear part
$$
a\cong a_xx+a_yy\quadand b\cong b_xx+b_yy\quad\text{modulo $(x,y)^2$}
$$ 
of our parameters.
Let  $c\in C$  be a rational point on the exceptional locus
for the blowing-up $q:T'\ra T$.
We assume that $c$ lies on the $a$-chart, the situation for the
$b$-chart being symmetric.
Using that $C_\red=\PP^1$, the ideal of $c\in T'$ is
generated by $x,y,z/a,b/a-\lambda$ for some scalar $\lambda\in k$.

\begin{proposition}
\mylabel{embedding dimension}
The  local ring $\O_{T',c}$ is regular if and only if the scalar $\lambda\in k$
satisfies the   equation $b_x+\lambda a_x+ b_y\lambda^2+a_y\lambda^3\neq 0$.
In any case,  $\edim(\O_{T',c})\leq 3$.
\end{proposition}

\proof
It follows from (\ref{a-chart}) that the embedding dimension
$\edim(\O_{T',c})$ is at most three, and that the $k$-vector space $\maxid_c/\maxid_c^2$
is generated by the classes of $x,y,b/a,z/a$ modulo the relations $\lambda x+y=0$ and $\lambda a+b=0$.
The latter relation equals $\lambda(a_xx+\lambda a_yy)+b_xx+\lambda b_y=0$ modulo $\maxid_c^2$. The assertion is now immediate.
\qed

\begin{corollary}
\mylabel{nonnormal}
The 2-dimensional scheme $T'$ is nonnormal if and only if both
parameters $a,b\in k[[x,y]]$ have no linear part.
\end{corollary}

\proof
We may assume that the ground field $k$ is separably closed, such that there
are infinitely many rational points on $C\simeq\PP^1$.
Then the scheme is $T'$ is nonnormal if and only if there are infinitely many $c\in C$ so that
$\O_{T',c}$ has embedding dimension three.
According to Proposition \ref{embedding dimension}, this holds if and only if $b_1+\lambda a_1+ b_2\lambda^2+a_3\lambda^3$
is the zero polynomial.  
\qed

\section{Maps to the Hilbert scheme}
\mylabel{maps hilbert}

We keep the notation from the preceding section, such that 
$q:S\ra T$ is the quotient morphism for an involution having
an isolated rational fixed point $s\in S$, with image $t\in T$.
Over the complement of the singularity $t\in T$, the quotient map
$q:S\ra T$ is a $G$-torsor, and in particular flat of degree two.
This gives an embedding $T\sm\left\{t\right\}\ra\Hilb^2(S)$ into the Hilbert scheme that
parameterizes   subschemes of length two. We may view this as a rational map $T\dra\Hilb^2(S)$.
Such rational maps extend to   morphisms on suitable blowing-ups of $T$. It turns out
that our blowing-up $g:T'\ra T$ defined in the preceding section already does the job.
To see this we have to come up with  a family of length two subschemes over $T'$,
which we do as follows:

Recall that $g:T'\ra T$ is the blowing-up of the primary ideal $(a,b,z)$
inside the local ring $\O_{T,t}^\wedge=k[[x,y,z]]$.
Now let $S'\ra S$ be the blowing-up of the induced ideal $(a,b,z)\O_{S,s}^\wedge=(a,b)\O_{S,s}^\wedge$
inside the local ring $\O_{S,s}^\wedge=k[[u,v]]$. 
The universal property of blowing-ups gives a morphism $h:S'\ra T'$, such that the   diagram
$$
\begin{CD}
S' @>>> S\\
@VhVV @VVqV\\
T' @>>g> T.
\end{CD}
$$
is commutative.

\begin{proposition}
The induced morphism $S'\ra T'\times S$ is a closed embedding,
and the projection $h:S'\ra T'$ is flat of degree two.
\end{proposition}

\proof
Let us first  check flatness:
The scheme $T'$ is integral outside the exceptional divisor for $T'\ra T$,
and clearly has no embedded components on the exceptional divisor $C\subset T'$.
Whence $T'$ is integral. The morphism $S'\ra T'$ is flat of degree two
over the complement of the exceptional divisor.
According to \cite{Hartshorne 1977}, Theorem 9.9, it is enough to
prove that  $h^{-1}(C_\red)\ra C_\red$ is flat of degree two.

It follows from (\ref{a-chart}) that  the $a$-chart
for the the exceptional divisor $C$ equals the spectrum of
$$
k[x,y,b/a,z/a]/(a,b,(z/a)^2+(b/a)^2x+y).
$$
Whence the reduction $C_\red$ has relations $x,y,z/a$, because
$a,b\in k[[x,y]]$ is a parameter system.
Consequently, the schematic preimage $h^{-1}(C_\red)$
is isomorphic to $\Spec k[u,v,b/a]/(u^2,v^2,z/a)$.
Using that $z/a=u\cdot b/a+v$, we see that the projection
$h^{-1}(C_\red)\ra C_\red$ is indeed flat of degree two.

Now let us check that $S'\ra T'\times S$ is a closed embedding.
This map is clearly proper, whence $\O_{S'}$ might be viewed as a coherent
sheaf on $T'\times S$. By the Nakayama Lemma, it 
suffices to show that $S'_{c}\ra S_c$
is a closed embedding for all points $c\in T'$. Making a field extension, we reduce to the case
that $c\in T'$ is rational.
There is nothing to prove if $c$ lies outside the exceptional divisor $C\subset T'$,
so let us assume $c\in C$. By symmetry, may also assume that $c$ lies in
the $a$-chart. Then there is a scalar $\lambda\in k$ so that the ideal defining $c\in C_\red$
is $(x,y,b/a-\lambda,z/a)\subset k[x,y,b/a,z/a]/(x,y,b/b,z/a)$.
The fiber $h^{-1}(c)\subset S'$  then is defined by
$$
(u^2,v^2,u\cdot b/a+v,b/a-\lambda)\subset k[u,v,b/a]/(u^2,v^2,u\cdot b/a+v,b/a-\lambda),
$$
which clearly defines a closed subscheme in $S$.
\qed

\medskip
It follows that our rational map $T\dra\Hilb^2(S)$ extends to a   morphism of schemes
$f:T'\ra\Hilb^2(S)$, which is defined by the family of subschemes $h:S'\ra T'$. The following fact
came as a surprise to me:

\begin{theorem}
The morphism $f:T'\ra\Hilb^2(S)$ is a closed embedding.
\end{theorem}

\proof
First observe that $f$ is proper: This is clear if $S$ is proper,
because then $T'$ is proper and the Hilbert scheme is separated. In general it follows  
by using a compactification $S\subset\overline{S}$.

Next, we check that the map $f:T'\ra\Hilb^2(S)$ is injective.
This is clear outside the exceptional divisor $C\subset T'$.
The $a$-chart of the reduction $C_\red$ is the spectrum of
$k[x,y,b/a,z/a]/(x,y,z/a)=k[b/a]$.
Given a closed point $c\in C$, say given by $b/a=\lambda$, the fiber in our family $h:S'\ra T'$
is $h^{-1}(c)=\Spec(k[u,v]/(u^2,u\lambda +v)$.
Clearly, different scalars $\lambda\in k$ give different subschemes $h^{-1}(c)\subset S$, whence
our map is indeed injective.

Let $h\in\Hilb^2(S)$ be a point with nonempty fiber $Y'_h$.
Using the Nakayama Lemma, it suffices to check that the fiber
$T'_h=f^{-1}(h)$  has length one. Making a base-change, we may assume that
the ground field is algebraically closed and that 
$h\in\Hilb^2(S)$ is closed.
Seeking a contradiction, we suppose that the fiber has length $>1$.
Then it contains a tangent vector $\Spec(k[\epsilon])\subset T'_h$,
where $\epsilon$ denotes an indeterminate subject to $\epsilon^2=0$.
Let $c\in T'$ be the support of such a tangent vector. Then the tangent
map $\Theta_{T'}(c)\ra\Theta_{\Hilb^2(S)}(h)$ is not injective.
We shall produce a contradiction by showing that the tangent map
actually is injective. We do this by producing a basis of $\Theta_{T'}(c)$ whose image
in $\Theta_{\Hilb^2(S)}(h)$ is linear independent.

To carry out this plan, let me recall the well-known description of the tangent space of 
the Hilbert scheme. For a nice discussion of these matters, see Artin's lecture notes (\cite{Artin 1976}, Section I.4)
or Vistoli's expository paper (\cite{Vistoli 1997}, Section 2).
Let $I\subset\O_{S}$ be the ideal of   $S_c\subset S$. Suppose $J\subset\O_S[\epsilon]$ is a coherent
ideal so that the quotient $\O_S[\epsilon]$ is $k[\epsilon]$-flat
and $J/\epsilon J=I$.
Choose generators $f_1=u^2$, $f_2=\lambda u+v$ for the ideal $I$.
Suppose $f_1',f_2'\in J$ are lifts for $f_1,f_2\in I$.
Then these lifts are necessarily generators of $J$. If $J'\subset\O_S[\epsilon]$ is another such ideal,
with lifts $f_1'',f_2''\in J'$, the
differences $f_i'-f_i''$ yields  an element in $\epsilon\O_S[\epsilon]=\epsilon\O_S$, and in turn
a residue class in $\epsilon\cdot\O_S/I$. It turns out that this gives
a well-defined homomorphism 
$$
\varphi:I/I^2\ra\epsilon\cdot\O_S/I,\quad f_i\longmapsto f_i'-f_i''.
$$
In this way, the tangent space $\Theta_{\Hilb^2(S)}(h)$ becomes
a torsor under the $k$-vector space $\Hom(I/I^2,\epsilon\cdot\O_S/I)$.
Throughout, we shall identify $\O_S$ as a subring of $\O_S[\epsilon]$,
and choose   $f''_1=f_1=u^2$, $f'_2=f_2=\lambda u+v$ as the obvious lifts.
This yields an identification of vector spaces
$$
\Theta_{\Hilb^2(S)}(h)\lra \Hom_k(I/I^2,\epsilon\cdot\O_S/I),\quad
\O_S/J\longmapsto \varphi.
$$

We are now ready for explicit computations:
Clearly our point $c$ contained in the exceptional locus $C\subset T'$ for the blowing-up $g:T'\ra T$.
By symmetry, we may assume that $c$ lies in the $a$-chart for the blowing-up.
Using the notation from  (\ref{a-chart}), we have $\maxid_c=(x,y,z/a,b/a-\lambda)$
for some scalar $\lambda\in k$. As explained in the proof for Proposition \ref{embedding dimension},
the $k$-vector space $\maxid_c/\maxid_c^2$ is generated by the residue classes of
$x,y,z/a,b/a-\lambda$, modulo the relations $x\lambda^2=y$ and $0=b+\lambda a$.

Consider the tangent vector $\psi:\maxid_c/\maxid_c^2\ra k$ that vanishes on the classes
of $x,y,b/z-\lambda$ and has  $\psi(z/a)=\epsilon$.
Then the fiber $S_\psi$ over $\Spec(k[\epsilon])\subset T'$
is isomorphic to the spectrum of $k[u,v,b/a]/(u^2,v^2,b/a-\lambda)$, and the
$k[\epsilon]$-algebra structure comes from $\epsilon \mapsto z/a$.
Using $z/a=b/a\cdot u+v$ and writing
$$
k[u,v,b/a]/(u^2,v^2,b/a-\lambda)=k[u,v,\epsilon]/(u^2,\epsilon +\lambda u+v),
$$
we see that   image of our tangent vector $f_*(\psi)\in\Theta_{\Hilb^2(A)}(h)$ corresponds to
the homomorphism with  $\varphi(u^2)=0$, $\varphi(\lambda u+v)=\epsilon$, which is obviously nonzero.
As a shorthand, we may represent the tangent vector $\psi$ as a $1\times 4$-matrix $(0,0,0,1)$ with respect
to the generating system $x,y,b/z-\lambda,z/a\in\maxid_c/\maxid_c^2$,
and its image $f_*(\psi)$ as the $1\times 2$-matrix $(0,\epsilon)$ with respect to the basis
$u^2,\lambda u+v\in I/I^2$.

Similar computations with other tangent vectors, which we leave to the reader,  yield the following data:
$$
\begin{array}[t]{c|c|c|c}
\psi\in\Hom(\maxid_c/\maxid_c^2,k)          &  (0,0,0,1)     &  (0,0,1,0)  &  (1,\lambda^2,0,0) \\
\hline&&&\\[-2ex]
 f_*(\psi)\in\Hom(I/I^2,\epsilon\cdot\O_S/I)  &  (0,\epsilon)  &  (0,\epsilon u)  &  (\epsilon, 0) \\
\end{array}
$$
Note that the last column is possible only if the embedding dimension of $T'$ is three, as explained
in Proposition \ref{embedding dimension}.
In any case, we see that the images $f_*(\psi)$ occurring in the second row are linearly independent.
The upshot is that the tangent map $f_*:\Theta_{T'}(c)\ra\Theta_{\Hilb^2(S)}(h)$ is injective.
\qed

\section{$G$-Hilbert schemes as partial resolutions}
\mylabel{g-hilbert schemes}

We keep the assumptions as in the preceding sections.
Let $\Hilb^G(S)\subset\Hilb^2(S)$ be the \emph{$G$-Hilbert scheme}. For me, this is the 
scheme that parameterizes $G$-invariant closed subschemes of length two on $S$.
Note that there are various other definitions in the literature, depending on the context at
hand.

The $G$-Hilbert scheme plays a central role in the \emph{McKay Correspondence}
for surface singularities over the complex numbers:
Ito and Nakamura \cite{Ito; Nakamura 1999} showed that the minimal resolution of 
rational double points is isomorphic to a suitable $G$-Hilbert scheme.
This was extended by Kidoh \cite{Kidoh 2001} to cyclic quotient singularities,
and by Ishii \cite{Ishii 2002} to arbitrary quotient singularities.
It turns out that the situation differs drastically in positive characteristics.
The goal of this section is to show that in our situation $\Hilb^G(S)$ usually contains
embedded components, and that the underlying reduced subscheme
is isomorphic to the blowing-up $T'\ra T$ constructed in
the preceding section.

Recall that the projection $S\ra T$ is a $G$-torsor over the complement
of the singularity $t\in T$.
It follows that the embedding $T\setminus\left\{t\right\}$ factors
over the $G$-Hilbert scheme.
Since $\Hilb^G(S)\subset\Hilb^2(S)$ is a closed embedding,
the closed embedding $T'\subset\Hilb^2(S)$ factors over
the $G$-Hilbert scheme as well.

\begin{proposition}
\mylabel{map bijective}
The closed embedding $T'\subset\Hilb^G(S)$ is bijective.
\end{proposition}

\proof
The subscheme $T\setminus\left\{t\right\}\subset\Hilb^G(S)$ parameterizes $G$-orbits
disjoint from the fixed point $s\in S$,
and its complement parameterizes $G$-invariant tangent vectors supported by $s\in S$.
It follows from Lemma \ref{fiber point} that any   tangent vector supported by $s\in S$ is $G$-invariant.
Consequently, the complement in question, viewed as a reduced subscheme, is isomorphic to the
projectivized cotangent space $\PP^1=\Proj(\Sym(\maxid_s/\maxid_s^2))$,
compare \cite{Brion; Kumar 2005}, Lemma 7.2.8.
By Theorem \ref{first blowing-up}, the reduced exceptional locus
of the blowing-up $T'\ra T$ is isomorphic to $\PP^1$ as well.
It follows that the closed embedding $T'\subset\Hilb^2(S)$ is bijective.
\qed

\medskip
Consequently, the reduction of the $G$-Hilbert scheme
is our  crepant partial resolution of singularities $T'=\Hilb^G_\red(S)$, which is possibly nonnormal, of the quotient surface $T=S/G$.
At this point I would like to point out that there is an a priori argument
that the normalization $\Hilb^G_\nor(S)$ is a partial resolution of singularities for the quotient surface $T$:
We may view the embedding $T\setminus\left\{t\right\}\subset\Hilb^G(S)$
as a rational map $T\dra\Hilb^G_\nor(S)$.
Suppose the inverse is undefined at some point of $\Hilb^G_\nor(S)$.
Then \cite{Beauville 1983a}, Lemma II.10 tells us that the rational map $T\dra\Hilb^G_\nor(S)$
contracts a curve. But we know that this rational map is an open embedding outside 
the closed point $t\in  T$, contradiction. 

\medskip
According to \cite{Brion; Kumar 2005}, Theorem 7.4.1 the Hilbert scheme $\Hilb^2(S)$ is smooth   and has
dimension four. It follows that the embedding dimensions of the $G$-Hilbert scheme 
are at most four.

\begin{theorem}
\mylabel{embedded component}
Suppose that both parameters $a,b\in k[[x,y]]$ have no linear terms.
Then $\Hilb^G(A)$ has an embedded component along the projective line
$\PP^1\subset\Hilb^G(A)$ that parameterizes tangent vectors  supported 
by the fixed point $s\in S$.
\end{theorem}

\proof
Let $h\in\Hilb^G(S)$ be a point corresponding to a tangent vector supported
by $s\in S$. We already know that $\Hilb^G_\red(S)=T'$ has embedding dimension three at $h$.
The idea now is to check that the local ring $\O_{\Hilb^G(S),h}$ has embedding dimension four.
Recall that the fixed scheme $S^G\subset S$ for the group action is defined by the
primary ideal $(a,b)\subset\O_{S,s}^\wedge$.  Obviously, we have closed embeddings
of Hilbert schemes $\Hilb^G(S)\supset\Hilb^2(S^G)\subset\Hilb^2(S)$,
and it follows from Lemma \ref{fiber point} that $h\in\Hilb^2(S^G)$.
In turn, we have  inclusions between tangent spaces 
$$
\Theta_{\Hilb^G(S)}(h)\supset\Theta_{\Hilb^2(S^G)}(h)\subset\Theta_{\Hilb^2(S)}(h).
$$
It suffices to prove that the inclusion on the right is bijective.
To this end, let $I\subset\O_{S,s}^\wedge$ be the ideal of the tangent vector
corresponding to $h\in\Hilb^2(S)$, and $\mathfrak{a}=(a,b)$. We have
\begin{gather*}
\Theta_{\Hilb^2(S)}(h)=\Hom(I/I^2,\O_S/I),\\
\Theta_{\Hilb^2(S^G)}(h)=\Hom((I+\mathfrak{a})/(I^2+\mathfrak{a}),\O_S/I).
\end{gather*}
Hence it suffices to check that the canonical surjection
$I/I^2\ra (I+\mathfrak{a})/(I^2+\mathfrak{a})$ is injective. The kernel
is $(I^2+\mathfrak{a})/I^2$, whence it remains to verify $\mathfrak{a}\subset I^2$.
We have inclusions
$$
(x,y)\subset (u^2,v^2)\subset(u,v)^2\subset I,
$$
where the first inclusion comes from Lemma \ref{fiber point}; 
for the last inclusion, see \cite{Brion; Kumar 2005}, Lemma 7.2.6.
Saying that $a,b\in k[[x,y]]$ have no linear parts means $(a,b)\subset(x,y)^2$,
which yields the desired inclusion $(a,b)\subset I^2$.
\qed

\section{Example: Rational double points}
\mylabel{example: rdp}

Keeping the same general assumptions as in the preceding sections,
we now consider the special case that the 2-dimensional singularity $\O_{T,t}$ is
a rational. For simplicity we also assume
that our ground field $k$ is algebraically closed.
By Corollary \ref{rational singularity},
the singularity is     a rational double point,
and our partial crepant resolution $T'=\Hilb^G_\red(S)$ is normal.

Recall that in characteristic zero, and for all primes $\geq 7$ as well,  rational double points
are \emph{taut} in the sense of Laufer (compare \cite{Laufer 1973a} and \cite{Laufer 1973b}).
This means that two rational double points are isomorphic if and only if the minimal
resolution of singularities have the same intersection graph, which in turn
correspond to the Dynkin diagrams.
As Artin computed in \cite{Artin 1977}, this does not hold true in characteristics
two, three, and five. For example, in characteristic two there are exactly $m$ different rational double
points   of Dynkin type $D_m$, which are denoted by $D_m^0,\ldots,D_m^{m-1}$, and five
different rational double points of Dynkin type $E_8$, which are denoted by $E_8^0,\ldots, E_8^4$.
Among other things, the isomorphism classes within a given Dynkin 
type differ by the \emph{Tyurina number}, which is length
of the scheme of nonsmoothness. Recall that if $f(x,y,z)=0$ defines an isolated singularity,
its   Tjurina number is the length of the \emph{Tjurina algebra}
$$
T=k[[x,y,z]]/(f,\frac{\partial f}{\partial x}, \frac{\partial f}{\partial y},\frac{\partial f}{\partial z}).
$$
Our first observation is that the amazing  diversity of rational double points in characteristic two 
does not occur among invariant rings for involutions:

\begin{proposition}
 \mylabel{d4r and e8}
In our situation, the rational double point $\O_{T,t}$ is  of type $D_{4r}^r$
for some integer $r\geq 1$, or of type $E_8^2$. The ring $\O_{T,t}^\wedge$ is isomorphic to  
 the quotient of the power series ring $k[[x,y,z]]$ by the polynomial $z^2+xy^rz+xy^{2r}+x^2y$ or by the polynomial
$z^2+yx^2z+x^5+y^3$, respectively.
\end{proposition}

\proof
In our situation, the local fundamental group $\pi_1^\loc(\O_{T,t})$ is cyclic of order two, and the universal covering
is given by the regular local ring $\O_{S,s}$.
On the other hand, Artin computed in \cite{Artin 1977}, Section 4  the local fundamental groups
of rational double points and gathered   information about some
unramified coverings as  well. We exploit his results as follows: First of all,
the singularities of type $A_m$ have no $2$-torsion at all in their local fundamental group,
which rules them out.
Among the rational double points of Dynkin type $E_m$, $m=6,7,8$, only $E_8^2$ has local fundamental
group of order two. This singularity is given, according to classification, by
the polynomial $z^2+yx^2z+x^5+y^3$. As we saw in Proposition \ref{normal form}, this equation indeed
describes the invariant ring of an involution on $\O_{S,s}=k[[u,v]]$.

It remains to treat the case of rational double points of Dynkin type $D_m$.
If $m=2n+1$ is odd, the singularity $D_m^r$ has  local fundamental group of order two
if and only if $n=2r$. However, 
the universal covering is then given by a rational double point of type $A_1$.
We conclude that $m=2n$ must be even. In case $r<m/4$, the singularity $D_{m}^r$ is simply connected.
In case $m/4<r$, the singularity $D_m^r$ admits a finite unramified covering by an $A_{8r-2m-1}$
singularity, which is never regular.
So the only remaining case is $m=4r$. 
By classification, rational double points of type $D_{4r}^r$ are given by the equation
$z^2+xy^rz+xy^{2r}+x^sy=0$. Indeed, their local fundamental group is
cyclic of order two. By Proposition \ref{normal form}, this equation actually
describes the invariant ring of an involution on $\O_{S,s}=k[[u,v]]$.
\qed

\medskip
Now let $T'=\Hilb^G_\red(S)$ be our crepant partial resolution of $T$ constructed in Section \ref{wild involutions},
which is the blowing-up whose center is the image of the fixed scheme $S^G\subset G$.
According to Corollary \ref{rational singularity}, the scheme $T'=\Hilb^G_\red(S)$ is normal, and the minimal resolution
$\tilde{T}\ra T$ factors over $T'$.
Whence $T'$ is obtained from $\tilde{T}$ by contracting all but one exceptional divisor.
We have to determine which exceptional divisor are contracted, and the isomorphism class
of  singularities  created on  $T'$. First, we treat the case that our rational double point 
$\O_{T,t}$ is of type $E_8$.

\vspace{2em}
\centerline{\includegraphics{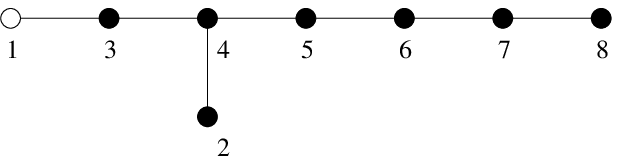}} 
\vspace{1em}
\centerline{Figure \stepcounter{figure}\arabic{figure}: The Dynkin diagram $E_8$.}
\vspace{1em}

The enumeration of vertices in the Dynkin diagram is always as in the Bourbaki tables \cite{LIE 4-6}.
These vertices corresponding to the exceptional curves $C_1,\ldots,C_8\subset\tilde{T}$.

\begin{proposition}
\mylabel{e8 singularity}
Suppose $\O_{T,t}$ is an $E_8^2$-singularity. Then
the partial resolution $T'\ra T$ is obtained from the minimal resolution $\tilde{T}\ra T$
by contracting all exceptional divisors except  $C_1\subset \tilde{T}$.
The singular locus of $T'$ consists of a rational double point of type $D_7^0$.
The situation is depicted in Figure 1.
\end{proposition}

\proof
We saw that $\O_{T,t}$ is given by the polynomial $z^2+yx^2z+x^5+y^3$, and $T'\ra T$ is the blowing-up
of the ideal $(y,x^2,z)$. We may decompose this blowing-ups into a blowing up of 
$(x,y,z)$, followed by a blowing-up of $(y/x,x,z/x)$ on the $x$-chart. 
Now recall the following simple but useful fact, which I learned from Torsten Ekedahl: 
For any rational double point, the blowing-up of
the reduced singular point introduces a single exceptional curve, and this curve
corresponds to the vertex in the Dynkin diagram adjacent to the \emph{longest root} of the root system
in question.
It follows that in the  blowing-up of an $E_8$-singularity the exceptional divisor corresponds
to $C_8$, and that the exceptional divisor for an iterated blowing-up corresponds to $C_1$,
confer the Bourbaki tables \cite{LIE 4-6}. We infer that $T'\ra T$ is obtained from the minimal resolution
by contracting all exceptional divisors except for $C_1$, and the singular locus of $T'$ 
consists of a rational double point of Dynkin type $D_7$.

To determine which $D_7^r$-singularity actually occurs, we use  Tjurina numbers.
We saw in the proof of Proposition \ref{first blowing-up} that the $a$-chart of $T'$ is generated by
$x,y,b/a,z/a$, modulo the two relations in (\ref{a-chart}). The scheme of nonsmoothness is thus defined by
the $2\times 2$-minors of the matrix of partial derivatives of the two relations.
A straightforward computations reveals that this determinantal scheme has length $12$.
According to \cite{Artin 1977}, page 15, the $D_7^r$-singularities have Tjurina number $12-2r$,
so our singularity must be of type $D_7^0$.
\qed

\medskip
Next, we consider the rational double points of type $D_{4r}$.

\vspace{2em}
\centerline{\includegraphics{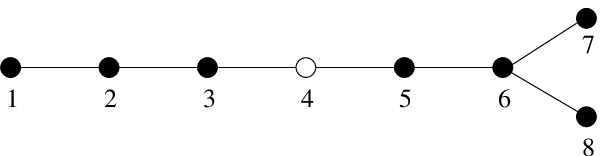} }
\vspace{1em}
\centerline{Figure \stepcounter{figure}\arabic{figure}: The   Dynkin diagram $D_8$.}
\vspace{1em}

\begin{proposition}
 \mylabel{d4r singularity}
Suppose $\O_{T,t}$ is a $D_{4r}^r$-singularity. Then the partial 
resolution $T'\ra T$  is obtained from the minimal resolution $\tilde{T}\ra T$
by contracting all exceptional divisors $C_i\subset\tilde{T}$, except the component $C_{2r}$.
The singular locus of   $T'$ consists of two rational double points, one  of type $A_{2r-1}$,
the other of type $D_{2r}^0$. In case $r=1$, we have to interpret $D_2^0$ as a pair of
$A_1$-singularities. The special case $r=2$ is depicted if Figure 2.
\end{proposition}

\proof
The singularity $\O_{T,t}$ is given by the equation
$z^2+xy^r+x^{2r}+yx^2$, and $T'\ra T$ is the blowing-up of
$(z,x,y^r)$. We may decompose this into a sequence  of $r$ blowing-ups or
reduced points, given on the $(i+1)$-th step by blowing-up the ideal $(z/y^i,x/y^i,y)$.
As in the preceding proof we infer that the last exceptional divisor
corresponds to the curve $C_{2r}\subset\tilde{T}$ on the minimal resolution.
Whence there are precisely two singularities on $T'$, one of type $A_{2r-1}$,
and one of type $D_{2r}$.
The Tjurina number on the $a$-chart of $T'$ turns out to be $2r$, which corresponds to
the $A_{2r-1}$-singularity. The Tjurina number on the $b$-chart is $4r$,
which implies that our singularity is of type $D_{2r}^0$. 
\qed

\section{Sign involution on abelian surfaces}
\mylabel{the sign involution}

In this section we apply our general results to a special case, which was my initial  motivation
to study Artin's wild involutions:
Let $A$ be an abelian surface over an algebraically closed field $k$ of characteristic two,
and $\iota:A\ra A$  be the sign involution. This gives an action of  the group 
$G=\sign$ on the abelian surface $A$.
The  fixed points  are precisely the 2-torsion points.
The kernel $A[2]\subset A$ of the multiplication-by-2 morphism $2:A\ra A$
is a group scheme of length $2^{2g}=16$, and necessarily of the form
$$
A[2] = (\ZZ/2\ZZ\oplus\mu_2)^\sigma\oplus N
$$
for some integer $\sigma=\sigma(A)$ subject to $0\leq\sigma\leq 2$,
and some local-local group scheme $N$ of length $4^{2-\sigma}$.
The integer $\sigma$ is called the \emph{$p$-rank} of the abelian surface,
and $A$ is called  \emph{ordinary} if $\sigma=2$. Another important invariant
is the embedding dimension   of $N$ , which is called the \emph{$a$-number} $a=a(A)$.
The $a$-number is also the dimension of  $\Hom(\alpha_p,A)$ viewed as a vector space
over $k=\End(\alpha_p)$, where $\alpha_p\subset\GG_a$ is the local additive group scheme of order $p$.
We have $0\leq a\leq 2$. Abelian surfaces is called \emph{supersingular} if  $a\geq 1$,
and \emph{superspecial} if $a=2$.

Shioda \cite{Shioda 1974} and Katsura \cite{Katsura 1978} studied the  singularities on the classical Kummer surface
$A/\left\{\pm 1\right\}$, in dependence  on $p$-rank and $a$-number.
The goal of this and the next section is to  complete  the  results of Shioda and Katsura
and determine the isomorphism class and   equations in normal
form for these singularities. Furthermore, we will determine the structure of the crepant partial resolution
 furnished by the $G$-Hilbert scheme.

Suppose first  that $A$ is not supersingular. Then
there are either two or four $2$-torsion points on $A$, and we have
to cope with  with the following slight complication: The $G$-Hilbert scheme $\Hilb^G(A)$
is no longer connected, because pairs of $2$-torsion point make up entire connected components.
However, the reduced connected component $\Hilb^{G,\circ}_\red(S)$ that is 2-dimensional
indeed yields a crepant partial resolution $T'=\Hilb^{G,\circ}_\red(S)$ of
the quotient surface $T=A/\sign$.

\begin{proposition}
\mylabel{singularity d4}
Suppose   $A$  is ordinary, that is, has $p$-rank $\sigma=2$.
Then the singular locus of
$A/\left\{\pm 1\right\}$ comprises four rational double points of type $D_4^1$.
Over each such singularity, the crepant partial resolution $\Hilb^{G,\circ}_\red(S)$ 
contains precisely three singularties, which are rational double points of type $A_1$.
\end{proposition}

\proof
The  singular points on $A/\left\{\pm 1\right\}$ are the images of
the $2$-torsion points on $A$, which are four in number.
According to Shioda \cite{Shioda 1974} and Katsura \cite{Katsura 1978}, Proposition 3,
each singularity is a rational double point  of type $D_4$.
These must be singularities of type $D_4^1$ by Proposition \ref{d4r and e8}.
The statement about the $G$-Hilbert scheme follows from Proposition \ref{d4r singularity}.
\qed

\medskip
The same arguments settle the following case as well:

\begin{proposition}
\mylabel{singularity d8}
Suppose the abelian surface $A$ has $p$-rank $\sigma=1$.
Then the singular locus of $A/\left\{\pm 1\right\}$ comprises two rational double points of type $D_8^2$.
Over each such singularity, the crepant partial resolution $\Hilb^{G,\circ}_\red(S)$ 
contains precisely two singularties, which are rational double points of type $A_3$ and $D_4^0$.
\end{proposition}

\medskip
The main task now is to understand the supersingular case, which is far more
challenging. Suppose $A$ is supersingular.
Then only the origin of $A$ is $2$-torsion, such that the normal surface
$A/\sign$ contains precisely one singularity.
Katsura showed in \cite{Katsura 1978},  Lemma 12 that the complete local ring at the singularity on $ A/\sign$ 
is isomorphic to $k[[x,y,z]]$ modulo the   polynomial  
\begin{gather*}
\label{katsura's equation}
f=q^4z^4+(1+(q^4-q)q^2x^3+q^2x^2y^2)z^2+ \\((q^4-q)x^3+q^2x^4y+x^2y^2)z +
(q^4-q)^2x^3+q^4x^5y^2+x^4y+xy^4,
\end{gather*}
for some parameter $q\in k$. This is, however, not  in Artin's normal form. Our first task
is to put Katsura's polynomial into Artin's normal form:

\begin{proposition}
\mylabel{supersingular equation}
Suppose $A$ is supersingular. Then the singularity on  $A/\left\{\pm 1\right\}$   is   formally isomorphic to the spectrum of
$k[[x,y,z]]/(z^2 +x^2bz + x^4y + xb^2)$,
where we have  $b=(q^4-q)x+y^2$ for some parameter  $q\in k$.
\end{proposition}

\proof
We simply  check that Katsura's polynomial $f$ is 
\emph{right equivalent} to our polynomial 
$$
g=z^2 +x^2((q^4-q)x+y^2)z + x^4y + x((q^4-q)x+y^2)^2.
$$
This simply means there is an automorphism of $k[[x,y,z]]$ sending $f$ to $g$.
Using   $\tilde{z}=z+q^2z^2$, we may rewrite Katsura's polynomial as
$$
f=\tilde{z}^2+ x^2((q^4-q)x+y^2)\tilde{z} + x^4y+x((q^4-q)x+y^2)^2+ x^4y(q^2z+q^4xy).
$$
Whence the inverse of the automorphism  $z\mapsto z+q^2z^2$ maps Katsura's polynomial $f$
to a power series of the form $g+x^4y\epsilon$ for some power series $\epsilon\in\maxid$.
So it remains to check that our polynomial $g$ and power series of the form $g+x^4y\epsilon$ are right equivalent.
One achieves this by inductively using substitutions of the form $y\mapsto y+y\epsilon$.
\qed

\medskip
The parameter $q\in k$ has the following geometric meaning:
Oort showed that any supersingular abelian surface $A$ is of the form
$(E\times E)/\alpha_2$, where $E$ is a supersingular elliptic curve,
and $\alpha_2\subset E\times E$ is an embedding of group schemes (\cite{Oort 1975}, Corollary 7).
Such embeddings depend on a single parameter $q\in\PP^1(k)$. If necessary, we may
interchange the factors in $E\times E$, and assume that $q\neq\infty$.
The resulting scalar $q\in k$ is precisely the parameter in our polynomial defining the
singularity. 

Let me now recall the following three facts: First, any product of $n\geq 2$ supersingular elliptic curves yields
isomorphic abelian varieties. In other words, there is only one
superspecial abelian variety in a given dimension $n\geq 2$. Second,   in characteristic two there is only one supersingular elliptic curve,
which is given by the Weierstrass equation $y^2=x^3+x$.
Third, $(E\times E)/\alpha_2$ is superspecial if and only if $q\in\FF_4$.

Having the equation for the singularity $T=A/\sign$, we now may easily infer the following
facts:

\begin{corollary}
\mylabel{supersingular normal}
Suppose the abelian surface $A$ is supersingular. Then the singularity on the classical Kummer surface
$A/\sign$ is not rational.
The crepant partial resolution $T'=\Hilb^G_\red(A)$ is normal if and only if $A$ is not superspecial.
\end{corollary}

\proof
Suppose $A$ is superspecial. Then the parameter $q\in k$ satisfies $q^4-q=0$,
and our equation defining the singularity reduces to $z^2+x^2y^2z+xy^4+yx^4$.
According to Corollary \ref{nonnormal}, the crepant partial resolution $T'$ is nonnormal.

Now suppose $A$ is not superspecial. Then $q^4-q\neq 0$, and we infer from Proposition \ref{embedding dimension}
that there is precisely one singularity on $T'$. But we saw in Section \ref{example: rdp}
that $T'$ would contain at least two singularities if the singularity on $A/\sign$ would be rational.
\qed

\medskip
To understand the situation, we have to enter   the theory of elliptic singularities,
which we do in the next section.

\section{Minimally elliptic singularities}
\mylabel{minimally elliptic}

Let me make a short digression and recall some   facts
on elliptic singularities. Suppose  $T=\Spec(\O_t)$ is a normal 2-dimensional local scheme,
with minimal resolution of singularities $f:\tilde{T}\ra T$,
and exceptional divisors $E_1\cup\ldots\cup E_n\subset \tilde{T}$.
The \emph{fundamental cycle} $Z=\sum n_iE_i$ is the smallest nonzero cycle supported
on the exceptional locus with integer coefficient so that $Z\cdot E_i\leq 0$.
The \emph{relative canonical cycle} $K=\sum d_iE_i$ is defined as the $\QQ$-valued cycle
supported on the exceptional locus satisfying the   system of linear equations
$Z\cdot E_i+E_i^2= -2\chi(\O_{E_i})$. Its coefficients $d_i\in \QQ$ are   called
the \emph{discrepancies}.
The next result due to Laufer (\cite{Laufer 1977}, Theorem 3.4 and
Theorem 3.10) is fundamental in the theory of surface singularities:

\begin{theorem}
\mylabel{properties elliptic}
The following are equivalent:
\renewcommand{\labelenumi}{(\roman{enumi})}
\begin{enumerate}
\item
$K=-Z$  as cycles holds.
\item 
The scheme $T$ is Gorenstein and the sheaf $R^1f_*\O_{\tilde{T}}$  has length one.
\item
We have $\chi(\O_Z)=0$, and $\chi(\O_{Z'})=-2$ for any subcycle $Z'\subsetneqq Z$.
\end{enumerate}
\end{theorem}

Singularities satisfying these equivalent condition are called \emph{minimally elliptic}.
They constitute a very interesting class of singularities, in importance
second only to   rational double points. They are special cases of \emph{elliptic singularities},
which are defined by the weaker condition $p_a=1$. Note, however, that  Reid  \cite{Reid 1997}
uses  slightly different terminology.

Now back to our main interest: Throughout this section we assume that 
$A$ is a supersingular abelian surface. 

\begin{proposition}
\mylabel{kummer elliptic}
The singularity on the classical Kummer surface $A/\sign$ is minimally elliptic.
If $A$ is not superspecial, then the singularity on the crepant partial  resolution
$T'=\Hilb^G_\red(A)$ is minimally elliptic as well.
\end{proposition}

\proof
Let $f:\tilde{T}\ra T$ be the minimal resolution of singularities of $T=A/\sign$.
By minimality, the relative canonical cycle satisfies $K_{\tilde{T}/T}\cdot E_i\geq 0$.
The singularity on $T$ is not rational, whence we actually have 
$K_{\tilde{T}/T}\cdot E_i>0$ for some exceptional divisor. We conclude that
$K_{\tilde{T}/T}<0$. Clearly $K_T=0$, whence $K_{\tilde{T}}<0$,
and therefore $H^2(\tilde{T},\O_{\tilde{T}})=0$.
This implies that the Picard scheme $\Pic_{\tilde{T}/k}$ is smooth, of expected dimension
$h^1(\O_{\tilde{T}})$. Let $\tilde{T}\ra P$ be the   Albanese morphism into the
abelian variety $P$ dual to $\Pic^0_{\tilde{T}/k}$.
Then the induced map $H^1(P,\O_P)\ra H^1(\tilde{T},\O_{\tilde{T}})$ is bijective.
In light of Proposition \ref{exceptional tree}, all exceptional curves $E_i\subset\tilde{T}$
are mapped to points on $P$. This means that  the boundary map $H^1(\tilde{T},\O_{\tilde{T}})\ra H^0(T,R^1f_*\O_{\tilde{T}})$
is zero. Now consider the exact sequence
$$
H^1(\tilde{T},\O_{\tilde{T}})\lra H^0(T,R^1f_*\O_{\tilde{T}})\lra H^2(T,\O_T).
$$
The term on the right is Serre dual to $H^0(T,\omega_T)$, which is 1-dimensional.
We conclude that $R^1f_*(\O_{\tilde{T}})$ is has length at most one.
It must have length one, because the   singularity on $T$ is nonrational.
The assertion on $T'=\Hilb^G_\red(A)$ follows in a similar way, because $\omega_{T'/T}$ is trivial.
\qed

\medskip
Now suppose that $A$ is not superspecial.
According to \cite{Katsura 1978}, the minimal 
resolution $\tilde{T}\ra T$ for the singularity on $T=A/\sign$ has the following intersection graph:

\vspace{1em}
\centerline{\includegraphics{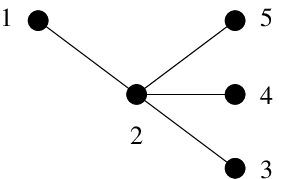}}
\vspace{.8em}
\centerline{Figure \stepcounter{figure}\arabic{figure}: Resolution graph for a minimally elliptic   double point.}
\vspace{1em}

\noindent
All components are isomorphic to $\PP^1$, and the  selfintersection numbers   are
$$
C_1^2=-3\quadand C_2^2=\ldots=C_5^2=-2.
$$
From this one directly computes  the fundamental cycle  
\begin{equation}
\label{fundamental cycle}
Z=-K=C_1+2C_2+C_3+C_4+C_5.
\end{equation}
The singularity appears in Laufer's Table 1 (\cite{Laufer 1977}, page 1288) under the name $A_{1,****}$.
It also appears in Wagreich's list of elliptic double points (\cite{Wagreich 1970}, Theorem 3.8)
under the symbol $\kreis{$\scriptstyle{4}$}^1_{0,1}$.
The singularity has multiplicity two,  according to \cite{Laufer 1977}, Theorem 3.13, because $Z^2=-1$.

By Proposition \ref{embedding dimension}, our crepant partial resolution $T'=\Hilb^G_\red(A)$
contains precisely one singularity $t'\in T'$. To describe it, we define an iterated
blowing-up of reduced points on $\tilde{T}$ as follows: First, blow-up a 
point on $C_2\subset\tilde{T}$ not contained in any other curves $C_i$.
Second, blow-up a point on the resulting exceptional divisor not contained in strict
transforms of any $C_i$. Call this two-fold blowing-up $\hat{T}\ra\tilde{T}$.
Let $\hat{C}_1,\ldots,\hat{C}_5\subset\hat{T}$ be the strict transforms of the $C_i\subset\tilde{T}$,
and let $\hat{C}_6,\hat{C}_7\subset\hat{T}$ be the two new exceptional curves.
They form the following intersection graph:

\vspace{1em}
\centerline{\includegraphics{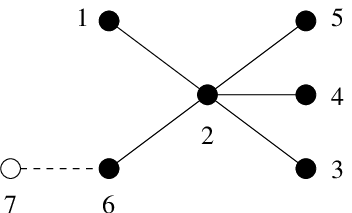}}
\vspace{.8em}
\centerline{Figure \stepcounter{figure}\arabic{figure}: Resolution graph for a minimally elliptic   triple point.}
\vspace{1em}

\noindent
The resulting selfintersection numbers are 
$$
\hat{C}_1^2=\hat{C}_2^2=-3,\quad \hat{C}_3^2=\ldots=\hat{C_6}^2=-2,\quadand \hat{C}_7^2=-1.
$$
The configuration of curve $\hat{C}_1\cup\ldots\cup\hat{C}_6\subset \hat{T}$ is negative definite
and contractible. (Let me remark in passing  that its contractibility would be a problem in characteristic zero,
compare \cite{Schroeer 2000}.)
 The resulting normal singularity is  minimally elliptic,
with fundamental cycle $\hat{Z}=\hat{C}_1+2\hat{C}_2+\hat{C}_3+\hat{C}_4+\hat{C}_5$,
which has $\hat{Z}^2=-3$. Whence this minimally elliptic singularity has embedding dimension three
and multiplicity three, by Laufer's result \cite{Laufer 1977}, Theorem 3.13.
It appears in loc.\ cit.\  Table 3, page 1293 under the   name
$A_{*,o}+A_{*,o}+A_{*,o}+A_{*,o}+A_{*,o}$.  

\begin{theorem}
Suppose $A$ is supersingular, but not superspecial. 
Then the singularity on the partial crepant resolution $T'=\Hilb^G_\red(A)$
is obtained as described above: Make a two-fold blowing-up $\hat{T}\ra\tilde{T}$
and contract  $\hat{C}_1\cup\ldots\cup\hat{C_6}\subset\hat{T}$.
The exceptional curve for $T'\ra A/\sign$ is the image of the $(-1)$-curve $\hat{C}_7\subset\hat{T}$.
\end{theorem}

\proof
Let $\hat{T}\ra T'$ be the minimal resolution of singularities, such that we have a commutative diagram
$$
\begin{CD}
\hat{T} @>>> T'\\
@VVV @VVV\\
\tilde{T} @>>> T.
\end{CD}
$$
Recall that $K_{T'/T}=0$ by Theorem \ref{first blowing-up}. On the other hand,
$K_{\tilde{T}}=K_{\tilde{T}/T}$ is the inverse of the fundamental cycle for $\tilde{T}\ra T$,
according to Proposition \ref{kummer elliptic}. We infer that $\hat{T}\ra \tilde{T}$ is not
an isomorphism. Instead, it factors into a sequence
$$
\hat{T}=\hat{T}_{n+5}\lra \hat{T}_{n+4} \lra \ldots\lra \hat{T}_6\lra \hat{T}_5=\tilde{T}
$$
of $n\geq 1$ blowing-ups with reduced points $t_{i+1}\in \hat{T}_i$ as centers. 
We now define curves $\hat{C}_i\subset\hat{T}$ for $1\leq i\leq n+5$ as follows:
For $1\leq i\leq 5$, let $\hat{C}_i$ be the strict transform of the exceptional curve $C_i\subset\tilde{T}$.
For $6\leq i\leq n+5$, let $\hat{C}_i$ be the strict transform of the exceptional curve for
the blowing-up $\hat{T}_{i}\ra\hat{T}_{i-1}$.
For convenience, we also denote by $\hat{C}_i\subset\hat{T}_j$ the images of $\hat{C}_i\subset\hat{T}$;
this ambiguity should not cause any confusion.

We have $t_{i+1}\in \hat{C}_i\subset\hat{T}_i$ by minimality of $\hat{T}\ra T'$.
Since the exceptional curve for $T'\ra T$ is irreducible, the morphism $\tilde{T}\ra T'$ contracts
  $\hat{C}_1\cup\ldots\cup\hat{C}_{n+4}\subset\hat{T}$ but not  $\hat{C}_{n+5}$.
We now use the following observation: Since the singularity $t'\in T'$ is minimally elliptic and
$K_{T'/T}=0$,
the multiplicities of $K_{\hat{T}_i}$ at the exceptional divisors $\hat{C}_1,\ldots,\hat{C}_{n+4}$ are negative,
and is zero at $\hat{C}_{n+5}$.

Suppose now that the first center $t_6\in \tilde{T}$ is contained in $C_1\setminus C_2$.
Then $K_{\hat{T}_6}$ has multiplicity zero along $\hat{C}_6$, which means $n=1$ and $\hat{T}=\hat{T}_6$.
We now obtain a contradiction as follows: By Corollary \ref{numerically cartier},
the Weil divisor $2g^{-1}(t)_\red\subset T'$, which is the image of $2\hat{C}_6$, is not numerically Cartier.  On the other hand, the vector
$$
\begin{pmatrix}
-4 & 1 & 0 & 0 & 0\\
1  & -2 & 1 & 1 & 1\\
0  & 1 & -2 & 0 & 0\\
0 &  1 & 0 & -2 & 0\\
0 &  1 & 0 & 0 & -2
\end{pmatrix}^{-1}\cdot
\begin{pmatrix}
2\\
0\\
0\\
0\\
0
\end{pmatrix}=
-\begin{pmatrix}
1\\
2\\
1\\
1\\
1
\end{pmatrix}.
$$
in integer-valued. This implies  that $2g^{-1}(t)_\red\subset T'$ is numerically Cartier, contradiction.
The matrix whose inverse appears above on the left is the 
intersection matrix $(\hat{C}_i\cdot\hat{C}_j)_{1\leq i,j\leq 5}$, 
and the vector on the left comprises the
intersection numbers   $(\hat{C}_i\cdot 2\hat{C}_6)_{1\leq i\leq 5}$.
In a similar way one excludes the cases $t_6\in C_i\setminus C_2$ for $i\geq 3$.
The case $t_6\in C_i\cap C_2$, $i\neq 2$ is impossible as well, because the
fiber $g^{-1}(t)\subset T'$ has multiplicity two by Theorem \ref{first blowing-up}.

Hence we have $t_6\in C_2$, and $t_5$ is  not contained in any $C_i$, $i\neq 2$.
Again using that $g^{-1}(t)\subset T'$ has multiplicity two, we infer
that $t_7\in \hat{C_6}$ is not contained in any other strict transform.
The canonical class $K_{\hat{T}_7}$ has multiplicity zero in $\hat{C}_7$.
As discussed above, this implies $\hat{T}=\hat{T}_7$,
and the assertion follows.
\qed

\section{The superspecial case}
\mylabel{superspecial case}

We now assume that our abelian surface $A$ is superspecial, that
is, isomorphic to $E\times E$, where $E$ is the supersingular elliptic curve,
which has Weierstrass equation $y^2=x^3+x$.
The minimally elliptic singularity $t\in T=A/\sign$ is formally given by the
equation $z^2+x^2y^2z+xy^4+yx^4=0$, according to Proposition \ref{supersingular equation}.
Katsura showed in  \cite{Katsura 1978} that the minimal resolution of singularities has the following
intersection graph:
 
\vspace{1em}
\centerline{\includegraphics{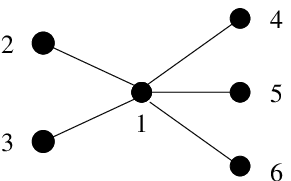}}
\vspace{.8em}
\centerline{Figure \stepcounter{figure}\arabic{figure}: Resolution graph for a minimally elliptic   double point.}
\vspace{1em}

\noindent
The intersection numbers are
$$
E_1^2=-3\quadand E_2^2=\ldots=E_5^2=-2.
$$
The fundamental cycle is $Z=2E_1+E_2+\ldots+E_5$, which has $Z^2=-2$.
Whence this minimally elliptic singularity has multiplicity two.
It appears in Laufer's Table 2 (\cite{Laufer 1977}, page 1290) under
the name $A_{*,o}+A_{*,o}+A_{*,o}+A_{*,o}+A_{*,o}$,
and in Wagreich's list (\cite{Wagreich 1970}, Theorem 3.8) under the symbol
$\kreis{$\scriptstyle{19}$}_0$.

Recall that our  crepant partial resolution $T'\ra T$
is given by the blowing-up of the primary ideal $(x^2,y^2,z)$.
We shall compare it with the blowing-up $T''\ra T$ of the maximal ideal
$(x,y,z)$. 

\begin{proposition}
\mylabel{trivial dualizing}
The schemes  $T',T''$ have trivial dualizing sheaf.
Both are nonnormal, and their normalizations are isomorphic.
The common normalization $S$ is obtained from the minimal resolution  of singularities $\tilde{T}$
by contracting all exceptional curves except $C_1\subset\tilde{T}$, compare Figure 5.
\end{proposition}

\proof
First, consider the blowing-up $T'\ra T$ of the primary ideal $(x^2,y^2,z)$.
Over the $x^2$-chart, the affine ring of $T'$ is given by  four generators $x,y,y^2/x^2,z/x^2$ modulo
two relations $(z/x^2)^2+x^2(y^2/x^2)(z/x^2) + x(y^2/x^2)^2 +y$ and
$(y^2/x^2)x^2=y^2$. The rational function $y/x$ is  clearly integral,
and the normalization  is given
by   three generators $x,y/x,z/x^2$ modulo the   relation
$(z/x^2)^2+x^2(y/x)^2(z/x^2)+ x(y/x)^4+x(y/x)$.

Now consider $T''\ra T$ be the blowing-up of the maximal ideal $(x,y,z)$.
Over the $x$-chart, the affine ring of  $T''$ is given  by three generators $x,y/x,z/x$ modulo
the relation $(z/x)^2 + x^3(y/x)^2(z/x) +x^3(y/x)^4 +x^3(y/x)$.
The rational function $z/x^2$ is clearly integral, and the normalization
is given by generators $x,y/x,z/x^2$ modulo the relation
$(z/x^2)^2+x^2(y/x)^2(z/x^2)+ x(y/x)^4+x(y/x)$. We conclude that
the schemes $T',T''$ have isomorphic normalizations $S$.
It is easy to check that $S$ contains five rational double points of type
$A_1$, whose resolution yield $\tilde{T}$.

It remains to verify the triviality of $\omega_{T''/T}$.
This can be done as in the proof for Theorem \ref{first blowing-up}.
\qed

\medskip
Next, we compute fibers over the singularity $t\in T$.
It is easy to see that $S_t$ is the trivial infinitesimal
extension of the projective line in the coordinate $y/x$ by the invertible
sheaf $\O_{\PP^1}(-2)$.  For simplicity, we   write $S_t=\PP^1\oplus\O_{\PP^1}(-2)$.

Next, recall that finite birational maps line $S\ra T'$
are determined by the conductor ideal $\mathfrak{c}\subset\O_S$,
which is the largest $\O_{T'}$-ideal that is at the same time an $\O_S$-ideal.

\begin{proposition}
\mylabel{conductor square}
The conductor ideals inside $\O_S$ for the normalizations
$S\ra T'$ and $S\ra T''$ coincide, and this conductor ideal is given by 
$\mathfrak{c}=\O_S(-S_t)\simeq\omega_S$.
\end{proposition}

\proof
First note $\omega_{T'}$ and $\omega_{T''}$ are trivial, whence
$\omega_{S/T'}=\omega_S=\omega_{S/T''}$.
By duality theory, the conductor ideals for $S\ra T'$ and $S\ra T'$ coincide
with the respective relative dualizing sheaves. We conclude that the two conductor ideals
coincide. To finish the proof, recall that the fundamental cycle of the minimally elliptic singularity $t\in T$
is $Z=2C_1+C_2+\ldots+C_6$, which implies   $\omega_S=\O_S(-S_t)$.
\qed

\medskip
To proceed, we have to compute the schematic image of the conductor scheme $S_t=\Spec(\O_S/\mathfrak{c})\subset S$:
We leave the following easy verification to the reader:

\begin{proposition}
\mylabel{conductor image}
The schematic image of $S_t=\PP^1\oplus\O_{\PP^1}(-2)$ on $T'$ is the fiber $T'_t=\PP^1\oplus\O_{\PP^1}(-1)$;
the induced morphism  between reduced subschemes is the Frobenius.
On the other hand, the schematic image of $S_t$ on $T''$ is the reduced fiber $(T''_t)_\red=\PP^1$;
the induced morphism between reduced subschemes is the identity.
\end{proposition}

\medskip
Summing up, the  two denormalizations $S\ra T'$ and $S\ra T''$ are given by the two cartesian and cocartesian squares
$$
\begin{CD}
\PP^1\oplus\O_{\PP^1}(-2) @>>> S\\
@Vg'VV @VVV\\
\PP^1\oplus\O_{\PP^1}(-1) @>>> T'
\end{CD}
\quad\quad\text{and}\quad\quad
\begin{CD}
\PP^1\oplus\O_{\PP^1}(-2) @>>> S\\
@Vg''VV @VVV\\
\PP^1 @>>> T''.
\end{CD}
$$
Here the glueing map $g'':\PP^1\oplus\O_{\PP^1}(-2)\ra\PP^1$ is just the identity on the underlying reduced
subschemes. In contrast, the glueing map $g':\PP^1\oplus\O_{\PP^1}(-2)\ra\PP^1\oplus\O_{\PP^1}(-1)$
is given by the relative Frobenius  morphism $\Fr:\PP^1\ra\PP^1$ on the underlying reduced subschemes.
Note that the relative Frobenius has  $\Fr^*\O_{\PP^1}(1)\simeq \O_{\PP^1}(2)$.
The passage from one denormalization $S\ra T''$ to another denormalization $S\ra T'$ 
might be called  an \emph{infinitesimal flip}.

We thus have completely unraveled the structure of $A/\sign$ and
its nonnormal crepant partial resolution $T'=\Hilb^G_\red(A)$.

\section{Serre conditions for symmetric products}
\mylabel{serre conditions}

In the forthcoming sections, we shall apply our results on Kummer surfaces to the geometry of Hilbert scheme
of points on abelian surface. To do this, it is first necessary to collect some facts
on symmetric products and their Cohen--Macaulay properties.
Fix a ground field $k$, for the moment of arbitrary characteristic $p\geq 0$,
and let $X$ be quasiprojective connected smooth scheme, say of dimension $g=\dim(X)$.
Throughout, we fix an integer  $n\geq 0$.
Recall that $n$-fold \emph{symmetric product} $\Sym^n(X)$ is the quotient of the $n$-fold product 
$X^n=X\times\ldots\times X$ by the action of the symmetric group
$S_n$ that permutes the factors. Such a quotient exists as a   scheme, and we have
$$
\Sym^n(X) = \bigcup \Spec(\TS^n(R)).
$$
Here the union runs over all affine open subsets $\Spec(R)\subset X$,
and $\TS^n(R)\subset T^n(R)$ is the subalgebra of symmetric tensors, which are by definition
the $S_n$-invariant tensors.
For an  account of symmetric tensors we refer to Bourbaki \cite{A 4-7}, Chapter IV, \S5.
The geometric points on $\Sym^n(X)$ correspond to formal linear combinations $\sum n_ix_i$ of pairwise different geometric points
with  $\sum n_i=n$, $n_i>0$.
As explained in Brion and Kumar's nice account (\cite{Brion; Kumar 2005}, Section 7.1), we have the following basic properties:

\begin{proposition}
\mylabel{basic properties}
The symmetric product $\Sym^n(X)$ is connected, quasiprojective, normal, $\QQ$-factorial,
of dimension $ng$, and its dualizing sheaf is invertible.
\end{proposition}

Throughout the paper, we follow the convention adopted by most researcher in the field
and call an algebraic scheme \emph{Gorenstein} if it is Cohen--Macaulay and
its dualizing sheaf   is invertible.
Note that our projective scheme $S=\Sym^n(X)$ is in general \emph{not} Cohen--Macaulay.
However, it   always has a dualizing sheaf, which for geometrically normal schemes might be defined
as the double dual $\omega_S=\det(\Omega^1_{S/k})^{\vee\vee}$.  For this dualizing sheaf, however, the Yoneda pairing
$$
\Ext^{ng-m}(\O_S,\shF)\times \Ext^m(\shF,\omega_S)\lra\Ext^{ng}(\O_S,\omega_S),
$$
together with the trace map, 
is known to give a perfect pairing only for $m=0$; compare \cite{Hartshorne 1977}, Chapter III, Section 7.

\begin{proposition}
\mylabel{gorenstein cohen--macaulay}
The symmetric product $\Sym^n(X)$ is Gorenstein if and only if it is
Cohen--Macaulay. These equivalent conditions hold provided
$p=0$ or $p>n$.
\end{proposition}

\proof
The first statement follows from the fact that the dualizing sheaf
$\omega_X$ is invertible, by Proposition \ref{basic properties}.
As to the second statement, suppose that $p=0$ or $p>n$. Then the order $n!$ of the symmetric group $S_n$
is invertible in the ground field $k$. According to Hochster and Eagon \cite{Hochster; Eagon 1971}, Proposition 13,
the quotient of the Cohen--Macaulay scheme $X^n$ by the $S_n$-action must again be Cohen--Macaulay.
\qed

\medskip
Note that   symmetric products are usually not Cohen--Macaulay in positive characteristics.
For example, Aramova \cite{Aramova 1992}, Proposition 2.8 computed explicitely that
$\Sym^n(X)$ for $g\geq 3$, $n\geq 2$ is not Cohen--Macaulay in characteristic two.
We now have a closed look into this matter.

Recall that a locally noetherian scheme
$S$ satisfies \emph{Serre's condition} $(S_k)$ if for each point $s\in S$,
the complete local ring $\O_{S,s}^\wedge$ is either Cohen--Macaulay, or contains a regular sequence
of length at least $k$, that is, $\depth(\O_{S,s}^\wedge)\geq k$.
The next result ensures that our symmetric products satisfy sufficiently many
Serre Conditions:

\begin{theorem}
\mylabel{symmetric serre}
The symmetric product $\Sym^n(X)$ satisfies Serre's Condition 
$(S_{g+2})$.
\end{theorem}

\proof
In the special cases $n\leq 1$ or $g\leq 1$, the symmetric product $\Sym^n(X)$ is smooth.
Hence it suffices to treat the case $n,g\geq 2$, such that $ng\geq g+2$.
Fix a point $s\in\Sym^n(X)$.
According to \cite{EGA IVb}, Corollary 6.7.2, we are allowed to extend our
ground field $k$, whence me may assume that $s$ is a rational point
of the form $s=\sum_{i=1}^r n_ix_i$, for certain rational points $x_i\in X$.
Our task now is to check that $\depth(\O_s^\wedge)\geq g+2$.

We first reduce with a standard argument to the case $r=1$:
Consider a preimage
$$
x=(\underbrace{x_1,\ldots,x_1}_{n_1},\underbrace{x_2,\ldots,x_2}_{n_2},\ldots, \underbrace{x_r,\ldots,x_r}_{n_r})\in X^n
$$
of $y\in\Sym^n(X)$.
As explained in \cite{Brion; Kumar 2005},  Lemma 7.1.3, the complete local ring
$\O_{s}^\wedge$ is isomorphic to the ring of invariants for the group  $G=S_{n_1}\times\ldots\times S_{n_r}$ 
inside $\O^\wedge_{x}$. It follows that
$$
\O_{s}^\wedge\simeq\O_{\Sym^{n_1}(X),n_1x_1}^\wedge\hat{\otimes}\ldots\hat{\otimes}\,\O_{\Sym^{n_r}(X),n_rx_r}^\wedge.
$$
This complete local ring is formally smooth if $n_1=\ldots=n_r=1$.
Hence it suffices to treat the case $n_1\geq 2$. By flatness of the tensor factors,
it suffices to check that the complete local ring of $n_1x_1\in\Sym^{n_1}(X)$ has
depth $\geq g+2$. In other words, we have reduced our problem to the case $r=1$.

We now suppose that our  point is of the form $s=nx_1$ for some rational point $x_1\in X$.
The corresponding point on $X^n$ is  $x=(x_1,\ldots,x_1)$.
Its complete local ring  is of the form
$$
\O_x^\wedge=
k[[\underbrace{u_{11},\ldots,u_{1g}}_g,\underbrace{u_{21},\ldots,u_{2g}}_g,\ldots,\underbrace{u_{n1},\ldots,u_{ng}}_g]],
$$
and the permutations $\sigma\in S_n$ act via $\sigma(u_{ij})=u_{\sigma(i),j}$.
This is also the completion of the symmetric algebra $\Sym(V)$ at the irrelevant ideal,
where $V=\oplus_{j=1}^gW$ is the $g$-fold sum of the
standard permutation representation of $S_n$ on $W=k^{\oplus n}$.
Clearly, the invariant subspace $W^{S_n}$ is 1-dimensional.
If follows that $V^{S_n}=\oplus_{j=1}^g(W^{S_n})$ has dimension $g$.
By the work of Ellingsrud and Skjelbred (\cite{Ellingsrud; Skjelbred 1980}, Theorem 3.9),
it follows that the irrelevant ideal of $\Sym(V)$ has depth $\geq g+2$.
\qed

\medskip
In dimension two, this tells us the following:

\begin{corollary}
Suppose $X$ is $2$-dimensional. Then $\Sym^2(X)$ is Cohen--Macaulay.
\end{corollary}

\proof
The $4$-dimensional scheme $\Sym^2(X)$ satisfies Serre's Condition $(S_4)$ by the preceding theorem,
whence is Cohen--Macaulay.
\qed

\medskip
In characteristic two, we may determine precisely what Serre Conditions hold. This also shows
that the preceding theorem gives the best general bound possible:

\begin{proposition}
\mylabel{2-fold symmetric}
Suppose $p=2$ and $g\geq 3$. Then the $2g$-dimensional scheme $\Sym^2(X)$ does not satisfy Serre's condition  $(S_{g+3})$.
\end{proposition}

\proof
We may assume that the ground field $k$ is algebraically closed.
Clearly, $\Sym^2(X)$ is Cohen--Macaulay outside the image of the fixed points
of the $S_2$-action on $X^2$. The fixed points are the diagonal points $(x_1,x_1)\in X^2$.
The complete local ring at a closed fixed point $x=(x_1,x_1)$ is of the form $\O_x^\wedge=k[[u_1,\ldots,u_g,v_1,\ldots,v_g]]$,
and $S_2$ acts via the involution $u_i\mapsto v_i$, $v_i\mapsto u_i$. This complete local ring is also 
the completion of the symmetric algebra $\Sym(V)$ at the irrelevant ideal,
with $V=k^{\oplus 2g}$, and group action given by the  block matrix
$$
\begin{pmatrix}
0 & \id_g\\
\id_g & 0
\end{pmatrix}
\in\GL(2g,k).
$$
Whence $V$ decomposes into $g$ irreducible 2-dimensional representations of $S_2=\ZZ/2\ZZ$.
According to Ellingsrud and Skjelbred \cite{Ellingsrud; Skjelbred 1980}, Corollary 3.2,
the depth of the 2g-dimensional invariant ring $\Sym(V)^{\ZZ/2\ZZ}$ at the irrelevant ideal
equals two plus the number of irreducible representations in $V$. The assertion follows.
\qed

\medskip
We also have the following negative result for higher dimensional schemes, which works in all characteristics:

\begin{proposition}
\mylabel{characterization cohen-macaulay}
Suppose $g\geq 3$. Then
the symmetric product $\Sym^n(X)$ is Cohen--Macaulay if and only if  
$n<p$.
\end{proposition}

\proof
As we already observed in Proposition \ref{basic properties}, 
the condition is sufficient by Hochster and Eagon \cite{Hochster; Eagon 1971}, Proposition 13.
For the converse, suppose that $n\geq p$. To proceed we may assume that the ground field
is algebraically closed.
Set $S=\Sym^n(X)$, and fix a point $s\in S$ of the form $s=nx_1$ for some closed point $x_1\in X$.
The complete local ring at the point $x=(x_1,\ldots, x_1)\in X^n$ is of the form
$$
\O_x=\hat{\bigotimes_{i=1,\ldots, n}}k[[x_1,\ldots,x_g]],
$$
 and the symmetric group $S_n$ acts via
permutations of the tensor factors. Let $P\subset S_n$ be a Sylow $p$-subgroup.
According to a sly computation of Campbell, Geramita, Hughes, Shank, and Wehlau 
(\cite{Campbell; Geramita; Hughes; Shank; Wehlau 1999}, Theorem 1.2), the invariant subring
$(\O_x^\wedge)^P\subset \O_x^\wedge$ cannot be Cohen--Macaulay.
This implies that the full invariant subring $(\O_x^\wedge)^{S_n}$ is not Cohen--Macaulay either.
The latter follows from the existence of a relative trace map $(\O_x^\wedge)^P\ra (\O_x^\wedge)^{S_n} $.
\qed

\medskip
Drawing from the huge amount of literature on depths in invariant rings, we may
generalize the preceding result  as follows:

\begin{proposition}
\mylabel{p-fold symmetric}
Suppose $\max(3,p)\leq n<2p$.  Then $\Sym^n(A)$ satisfies Serre's condition
$(S_{g+2})$, but not $(S_{g+3})$.
\end{proposition}

\proof
We may assume that the ground field $k$ is algebraically closed.
Fix a closed point $y=\sum_{i=1}^r n_ia_i$ in $S=\Sym^n(A)$.
We first consider the case  
$r=1$. Then $y$ is the image of $x=(a,\ldots,a)\in A^n$ with $a=a_1$.
As in the preceding proof, the complete local ring $\O^\wedge_{S,y}$
is the formal completion of the symmetric algebra $\Sym(V)$, where $V=(k^{\oplus g})^{\otimes n}$,
and the symmetric group $S_n$ acts via permutations of the tensor factors.
We now invoke a result of Kemper (\cite{Kemper 2001}, Theorem 3.3)  and have to verify some
hypothesis. First of all, our condition   $p\leq n<2p$ ensures that 
the Sylow $p$-subgroup $P\subset S_n$ has order $p$,
and $P$ equals its own normalizer. This $S_n$-action on $V$ 
leaves the standard basis  invariant, and the Sylow $p$-subgroup $P$
decomposes this standard basis into $g$ orbits.
According to loc.\ cit., the depths of the invariant
subring $\Sym(V)^{S_n}$ at the irrelevant ideal is $\min(g+2,ng)$.
In light of $3\leq n$, this minimum equals $g+2$.
We infer that $\depth(\O_{S,s})=g+2$. It follows that the local ring
$\O_{S,s}$ satisfies $(S_{g+2})$ but not $(S_{g+3})$.

It remains to treat the general case $r\geq 1$.
As explained in \cite{Brion; Kumar 2005},  Lemma 7.1.3, the complete local ring
$\O_{S,s}^\wedge$ is isomorphic to the ring of $S_{n_1}\times\ldots\times S_{n_r}$-invariants
inside $\O^\wedge_{X^n,(x_1,\ldots,x_n)}$. It follows that
$$
\O_{S,s}^\wedge\simeq\O_{\Sym^{n_1}(A),n_1a_1}^\wedge\hat{\otimes}\ldots\hat{\otimes}\,\O_{\Sym^{n_r}(A),n_ra_r}^\wedge.
$$
If all factors have $n_i<p$, then the product  is Cohen--Macaulay, whence its depth equals $ng\geq g+2$.
If some factor has $p\leq n_i$, we are in case $r=1$, so we also have   depth $\geq g+2$.
In any case, the the local ring $\O_{S,s}$ satisfies Serre's condition $(S_{g+2})$.
\qed

\medskip
I close this section  with an observation concerning  rational  singularities.
Let $Y$ be a normal irreducible scheme of finite type that admits a resolution
of singularities $f:X\ra Y$.
Recall that $Y$ is said to have only  \emph{rational singualrities}
if  $R^if_*(\O_X)=0$ and $R^if_*(\omega_X)=0$ for all $i>0$. 
Note that the condition on the higher direct images of the dualizing sheaf
are superfluous in characteristic zero, thanks to the Grauert--Riemenschneider Vanishing Theorem
\cite{Grauert; Riemenschneider 1970}.

Now suppose that $Y$ has only rational singularities. Then
the shifted sheaf $f_*(\omega_X)[d]$, $d=\dim(Y)$ is a dualizing complex on $Y$, and in particular
$Y$ must be Cohen--Macaulay. 
As a consequence of the preceding two propositions, we obtain the following fact:

\begin{corollary}
\mylabel{symmetric nonrational}
Suppose either $g\geq 3$ and $n\geq p$, or $\max(3,p)\leq n<2p$. Then the  symmetric product $\Sym^n(X)$  are has not only rational singularities.
\end{corollary}

\section{Symmetric products of abelian varieties}
\mylabel{symmetric products}

Now let $A$ be a $g$-dimensional abelian variety. In this section we study
its symmetric products $\Sym^n(A)$.
The first thing to say is that, since $\omega_A$ is trivial,
the dualizing sheaf of $\Sym^n(A)$ is trivial as well,
as explained in \cite{Brion; Kumar 2005}, Section 7.1.
Another special feature of the situation is the \emph{addition map}
$$
+:A^n\lra A,\quad (a_1,\ldots,a_n)\longmapsto a_1+\ldots+a_n.
$$
It is clearly $S_n$-invariant, whence descends to an addition map $\Sym^n(A)\ra A$.
Given a point $t\in A$, we denote by $\Sym^n_t(A)$ the fiber over $t\in A$ of
the addition map. 
 More generally, if $\phi:T\ra A$ is a morphism of schemes, we define
$\Sym^n_\phi(A)$ via the pull back
$$
\begin{CD}
\Sym^n_\phi(A) @>>>\Sym^n(A)\\
@VVV @VV+V\\
T @>>\phi>  A.
\end{CD}
$$
We shall be particularly interested in the  pullback under the multiplication-by-$n$ map $n:A\ra A$.
Therefore, I formulate the next result in the  abstract language of $T$-valued points:

\begin{proposition}
\mylabel{pullback isomorphic}
Let $T$ be a scheme, and $\phi,\psi:T\ra A$ be two morphisms.
Suppose that the group element $\phi-\psi\in A(T)$ is $n$-divisible.
Then the choice of $d\in A(T)$ with $\phi-\psi=nd$ induces an
isomorphism of $T$-schemes $\Sym_\psi^n(A)\ra\Sym_\phi^n(A)$.
\end{proposition}

\proof
The morphism $d:T\ra A$ corresponds to the  section $d\times\id$ for the projection
$A\times T\ra T$. It yields a translation  
$$
\tau_{d}:A\times T\lra A\times T,\quad
(a,t)\longmapsto (a+d(t),t),
$$
which is a $T$-morphism.
This induces a $T$-morphism $\Sym^n(\tau_{d})$, such that the diagram
$$
\begin{CD}
\Sym^n(A)\times T @>\Sym^n(\tau_d)>>\Sym^n(A)\times T\\
@V+VV @VV+V\\
A\times T @>>\tau_{nd}> A\times T
\end{CD}
$$
is commutative.
Similarly, $\phi,\psi:T\ra A$ yield sections for the projection $A\times T\ra T$,
and we have a commutative diagram
$$
\begin{CD}
A\times T @>\tau_{nd}>> A\times T\\
@A\psi\times\id AA @AA\phi\times\id A\\
T @>>\id_T> T.
\end{CD}
$$
Hence $\id_T\times\Sym^n(\tau_{d})$ induces an isomorphism 
$$
(A,\psi\times\id)\times_{(A\times T)}(\Sym^n(A)\times T)\lra
(A,\phi\times\id)\times_{(A\times T)}(\Sym^n(A)\times T).
$$
Identifying these fiber products with $\Sym^n_\psi(A)$ and $\Sym^n_\phi(A)$, respectively,
we obtain the desired result.
\qed

\medskip
Applying this with the multiplication-by-$n$ map, we reach the following:

\begin{corollary}
\mylabel{pullback product}
Let pull-back $\Sym^n(A)\times_A (A,n)$ with respect to the multiplication-by-$n$ map $n:A\ra A$
is canonically isomorphic to the product $\Sym^n_0(A)\times_k A$.
\end{corollary}

\proof
Set $T=A$. Let $\psi$ be the multiplication-by-$n$ map $n:A\ra A$,
and   $\phi$ be the zero map $0:A\ra A$.
Then $\psi-\phi=n\id_A$ is obviously $n$-divisible, so we may apply Proposition \ref{pullback isomorphic}.
\qed

\medskip
This means  that the $A$-scheme $\Sym^n(A)$ is a \emph{twisted form}
of the  product $A$-scheme $\Sym^n_0(A)\times_k A$, with respect to the finite flat topology.
In particular, the schemes $\Sym^n_\eta(A)$ over the function field $\kappa(A)$ is a twisted form
of $\Sym^n_0(A)\otimes_k\kappa(A)$.
Note that singularities may dissappear upon passing to twisted forms, as explained in 
\cite{Schroeer 2006b}.

Next, let $A^n_0\subset A$ be the kernel of the addition map $A^n\ra A$,
which is an abelian variety of dimension $ng-1$.
Clearly, this kernel is invariant under the permutation action of $S_n$.
So we may form the quotient scheme $A^n_0/S_n$, which is normal,
and obtain a morphism $A^n_0/S_n\ra\Sym^n(A)$, which factors over the closed fiber
$\Sym^n_0(A)\subset\Sym^n(A)$ of the addition map. 

\begin{proposition}
\mylabel{base-change}
Suppose $g\geq 2$. Then the canonical morphism $A^n_0/S_n\ra \Sym^n_0(A)$ is an isomorphism.
\end{proposition}

\proof
We may assume that $k$ is algebraically closed.
Let $U\subset A^n$ be the open subset whose closed points are
the $(a_1,\ldots,a_n)$ with pairwise different entries, and let
$V\subset\Sym^n(A)$ be its image. 
Then for all closed points $u\in U$, the stabilizer $G_u\subset G$ is
trivial. Then the  projection $U\ra V$ is flat of degree $n!$, and   the formation of the quotient
$V=U/S_n$ commutes with arbitrary base change in $V$. Using that $U_0=A^n_0\cap U$ is the preimage
of $V_0=\Sym^n_0(A)\cap V$, we deduce that the morphism
$U_0/S_n\ra V_0$ is an isomorphism.
Clearly, the complement $A^n_0\setminus U\subset A^n_0$ has codimension $g$, and we have $g\geq 2$
by assumption.
We infer that the morphism $A^n_0/S_n\ra\Sym^n_0(A)$ is an isomorphism in codimension $\leq 1$.

Since $A^n_0/S_n$ is normal, we see that $\Sym^n_0(A)$ is regular in codimension $\leq 1$.
Moreover,   $\Sym^n_0(X)$ satisfies Serre's Condition $(S_2)$
by Theorem \ref{symmetric serre}. 
According to Zariski's Main Theorem, the finite morphism   $A^n_0/S_n\ra \Sym^n_0(A)$
must be an isomorphism.
\qed

\medskip
As an application, we infer that the generic fiber $\Sym^n_\eta(A)$   for the addition map contains
no hidden singularities. Note that the corresponding statement for Hilbert schemes does not hold,
as we shall see  in Section \ref{hilbert-chow}.

\begin{corollary}
\mylabel{regular smooth}
Suppose $g\geq 2$, and let $y\in\Sym^n_\eta(A)$ be a point.
Then the local ring $\O_{y}$ is regular if and only if it is geometrically regular
as $\kappa(\eta)$-algebra.
\end{corollary}

\proof
The condition is obviously sufficient.
To check that it is also necessary, set $F=\kappa(\eta)$ and choose
an algebraic closure $F\subset\bar{F}$.
According to Corollary \ref{pullback product}, the geometric generic fiber $\Sym^n_\eta(A)\otimes_F \bar{F}$
is isomorphic to $\Sym^n_0(A)\otimes_k \bar{F}$, and the latter
is isomorphic to $A^n_0/S_n\otimes\bar{F}$ by the preceding theorem.
Set $B=A^n_0$ and $G=S_n$, and let $y\in B/G$ a point whose local ring $\O_{B/G,y}$
is regular. Since $B$ is Cohen--Macaulay, the projection $f:B\ra B/G$ is flat near $y$, whence
the schematic fiber $f^{-1}(y)$ has length $n$.
Then the induced morphism $B\otimes\bar{F}\ra B/G\otimes\bar{F}$ is flat at the preimages of $y$.
Since $B\otimes\bar{F}$ is regular, the scheme $B/G\otimes\bar{F}$ must be regular at the preimages of $y$.
\qed

\medskip
We now turn to the special case $n=2$.
Then the antidiagonal $A\ra A^2$, $a\mapsto (a,-a)$ yields an isomorphism $A\ra A^2_0$,
and the permutation action of the symmetric group $S_2$ restricts to the action of
$\left\{\pm 1\right\}$ via the sign involution
$a\mapsto -a$.
Summing up, we have:

\begin{corollary}
\mylabel{involution symmetric}
For $g\geq 2$ we have an isomorphism $A/\left\{\pm 1\right\}\ra\Sym^2_0(A)$.
\end{corollary}

This allows to apply our results on the sign involution from Section \ref{the sign involution}
to twofold symmetric products.

\section{The Hilbert--Chow morphism and quasifibrations}
\mylabel{hilbert-chow}

Let $k$ be a ground field of arbitrary characteristic $p\geq 0$.
Recall that for any $k$-scheme of finite type $X$, the Hilbert scheme
is related to   symmetric products via the
\emph{Hilbert--Chow} morphism
$$
\gamma:\Hilb^n(X)\ra\Sym^n(X).
$$
It sends a finite subscheme $A\subset X$ to the   sum of points
$\sum n_ix_i$, where the coefficients are the lengths of the Artin rings $\O_{A,x_i}$.
In my opinion, the fact that such a map exists as a morphism of schemes is rather nontrivial.
Iversen \cite{Iversen 1970} worked this out, using his theory of linear determinants.
If $X$ is a smooth connected surface, then $\Hilb^n(X)$ is again smooth and connected, of dimension $2n$, and
the Hilbert--Chow morphism
$\gamma:\Hilb^n(X)\ra\Sym^n(S)$ is a crepant resolution of singularities,
as explained in  \cite{Brion; Kumar 2005}, Section 7.4.

Now let $A$ be an abelian surface.
Then $\Hilb^n(A)$ is smooth and connected, and its dualizing sheaf is trivial.
We may compose the Hilbert--Chow morphism
with the addition map and obtain another addition map
$$
\Hilb^n(A)\stackrel{\gamma}{\lra}\Sym^n(A)\stackrel{+}{\lra} A.
$$
As with symmetric products, we denote by $\Hilb^n_t(A)$, $t\in A$ its fibers. 
More generally, if $\psi:T\ra A$ is a morphism,
we define $\Hilb^n_\psi(A)$ as the corresponding base change.
The analogue of Proposition \ref{pullback isomorphic} holds true for Hilbert schemes,
with the same proof:

\begin{proposition}
\mylabel{pullback hilbert}
Let $T$ be a scheme, and $\phi,\psi:T\ra A$ be two morphisms.
Suppose that the group element $\phi-\psi\in A(T)$ is $n$-divisible.
Then the choice of $d\in A(T)$ with $\phi-\psi=nd$ induces an
isomorphism  $\Hilb_\psi^n(A)\ra\Hilb_\phi^n(A)$ of $T$-schemes.
\end{proposition}

As a  consequence:

\begin{corollary}
\mylabel{addition flat}
The addition map $f:\Hilb^n(A)\ra A$ is flat,
and the canonical map $\O_A\ra f_*(\O_{\Hilb^n(A)})$ is bijective. 
\end{corollary}

\proof
By Proposition \ref{pullback hilbert}, the pull-back $\Hilb^n(A)\times_A (A,n)$ with respect to the
multiplication-by-$n$ map $n:A\ra A$ is isomorphic to the product $\Hilb^n_0(A)\times_k A$.
This implies flatness of the addition map.
To proceed, consider the commutative diagram:
$$
\xymatrix{ 
 \Hilb^n(A) \ar[r]^\gamma \ar[dr]_f& \Sym^n(A) \ar[d]_g & A^n \ar_{\phantom{qqqq} q} [l]\ar[dl]^h \\
 & A  }
$$
Here the lower arrows are the addition maps.
The composition of the two inclusions $\O_A\subset g_*(\O_{\Sym^n(A)})\subset h_*(\O_{A^n})$ is bijective,
whence $\O_A=g_*(\O_{\Sym^n(A)})$. This implies that $\Sym^n_0(A)$ is connected.
As explained in \cite{Brion; Kumar 2005}, Section 5, the Hilbert--Chow morphism
$\gamma:\Hilb^n(A)\ra\Sym^n(A)$ is birational. Using that $\Sym^n(A)$ is normal, we infer
with Zariski's Main Theorem that $\O_{\Sym^n(A)}=\gamma_*(\O_{\Hilb^n(A)})$.
The result follows.
\qed 

\medskip
The splitting of $\Hilb^n(A)\times_A (A,n)$ is in line with 
Beauville' splitting result on K\"ahler manifolds
with zero Ricci curvature (\cite{Beauville 1983}, Theorem 1).
Let us now take a closer look at the fiber over the origin, which is
Beauville's  \emph{generalized Kummer variety} $\Km^n(A)=\Hilb^n_0(A)$:

\begin{proposition}
\mylabel{integral intersection}
Beauville's generalized Kummer variety $\Km^n(A)$ is integral, locally of complete intersection,
and has trivial dualizing sheaf.
The induced   morphism $\gamma_0:\Km^n(A)\ra\Sym^n_0(A)$
is birational, and a crepant partial   resolution.
\end{proposition}

\proof
The   fiber $\Hilb^n_\eta(A)$ over the generic point $\eta\in A$ is obviously integral, regular  
 hence locally of complete intersection, and has trivial dualizing sheaf.
Moreover, the morphism $\gamma_\eta:\Hilb^n_\eta(A)\ra\Sym^n_\eta(A)$ is a crepant resolution of singularities. 
Over the algebraic closure of the function field $\kappa(\eta)$,
this schemes and morphism become isomorphic to $\Km^n(A)$ and $\gamma_0:\Km(A)\ra\Sym^n_0(A)$,
after base-changing to $\overline{\kappa(\eta)}$.
The assertions now follow from descend theory.
\qed

\medskip
We already saw that the canonical map $A^n_0/S_n\ra\Sym^n_0(A)$ is an isomorphism.
In particular, $\Sym^n_0(A)$ is normal. In contrast, $\Km^n(A)$ is not necessarily normal,
as we shall see.
To understand this, let us consider the case $n=2$.
The group $G=\sign$ acts on the abelian surface $A$ via the sign involution.
Recall that  $\Hilb^{G,\circ}_\red(A)\subset\Hilb^2(A)$ is the 2-dimensional integral component
inside the fixed scheme $\Hilb^2(A)^G$.

\begin{proposition}
\mylabel{identification}
We have $\Hilb^{G,\circ}_\red(A)=\Km^2(A)$ as subschemes of $\Hilb^2(A)$.
\end{proposition}

\proof
Since both subschemes are integral, it suffices to check that both
are mapped onto $\Sym^2_0(A)$ under the Hilbert--Chow morphism $\gamma:\Hilb^2(A)\ra\Sym^2(A)$,
which is obvious.
\qed

\medskip
So far, the results hold true in arbitrary characteristics. We now turn to the
case of characteristic two:

\begin{theorem}
Suppose $p=2$. Then all fibers of the addition map
$\Hilb^2(A)\ra A$ are geometrically integral, but nonsmooth.
They are geometrically normal if and only if the abelian surface $A$
is not superspecial. In particular, this applies to Beauville's
generalized Kummer variety $\Km^2(A)=\Hilb^2_0(A)$.
\end{theorem}

\proof
It follows from Proposition \ref{pullback hilbert} that for any point $t\in A$,  say with residue field $F=\kappa(t)$,
there exists an isomorphism  $\Hilb^2_t(A)\otimes_F \bar{F}\simeq \Km^2(A)\otimes_k \bar{F}$.
Hence it suffices to treat the zero fiber $\Km^2(A)$.
By Proposition \ref{identification}, we have an identification of $\Km^2(A)$
with the $G$-Hilbert scheme $\Hilb^{G,\circ}_\red$. On the other hand,
by Proposition \ref{map bijective} we have  $T'=\Hilb^{G,\circ}_\red$,
where $T'$ is the canonical blowing-up attached to the quotient $T=A/\sign$.
We analyzed this partial crepant resolution in 
Sections \ref{the sign involution}, \ref{minimally elliptic}, and \ref{superspecial case}.
In all cases, $T'$ is     integral, but not regular. 
It is normal if and only if $A$ is not superspecial.
\qed

\medskip
Let me call a morphism $f:X\ra Y$ between smooth proper connected schemes
with $\O_Y=f_*(\O_X)$ a \emph{quasifibration} if the generic fiber
$X_\eta$ is not smooth. The generic fiber $X_\eta$ is always a regular scheme.
In characteristic zero, this implies that $X_\eta$ is smooth over $\kappa(\eta)$,
so there are no quasifibrations. However, there are quasifibrations in positive
characteristics. The most prominent are the quasielliptic surfaces in characteristic two and three.
We just saw that the addition map $\Hilb^2(A)\ra A$ is a quasifibration in characteristic two.
If $A$ is supersingular but not superspecial, the geometric generic fiber contains
a minimally elliptic singularity. Recently,  Hirokado \cite{Hirokado 2004} studied
quasifibrations with simple elliptic singuarities.
It would be interesting to determine under what conditions on $p$, $n$, and $A$
the addition map $\Hilb^n(A)\ra A$ is a quasifibration.

Let me close this paper with an observation on canonical singularities:
Suppose $Y$ is a connected normal scheme of finite type, admitting a 
resolution of singularities $f:X\ra Y$. 
One says that $Y$ has only \emph{canonical singularities} $Y$ if the reflexive rank-one sheaf
$\omega_Y^{[r]}=(\omega^{\otimes r})^{\vee\vee}$ is invertible for some integer $r\geq 1$,
and the canonical inclusion $f_*(\omega_X^{\otimes r})\subset\omega_Y^{[r]}$ is bijective.
Elkik \cite{Elkik 1981} proved that  canonical singularities in characteristic zero 
are rational.
This does not hold true in positive characteristics. Indeed, our symmetric products
of abelian surfaces yield counterexamples:

\begin{proposition}
\mylabel{canonical singularities}
Suppose $\max(3,p)\leq n<2p$. Then the singularities of the  symmetric product $\Sym^n(A)$ are
canonical, but not rational. 
\end{proposition}

\proof
We saw in Proposition \ref{symmetric nonrational} that   $\Sym^n(A)$ has not only rational
singualrities.
On the other hand, this scheme is normal with $\omega_{\Sym^n(A)}=\O_{\Sym^n(A)}$,
and the Hilbert--Chow morphism $\gamma:\Hilb^n(A)\ra\Sym^n(A)$ is a crepant resolution.
Hence
$$
\gamma_*(\omega_{\Hilb^n(A)})=\gamma_*(\O_{\Hilb^n(A)})=\O_{\Sym^n(A)}=\omega_{\Sym^n(A)}
$$
trivially holds.
\qed


\end{document}